\documentclass{amsart}

\usepackage{fullpage,url,amssymb,enumerate,colonequals}
\usepackage{url,amssymb,enumerate,colonequals}
\usepackage{multirow}
\usepackage{mathrsfs} % for \mathscr (script letters)
\usepackage[section]{placeins}
\usepackage{MnSymbol}
\usepackage{extarrows}
\usepackage{lscape}
\usepackage[all,cmtip]{xy}

\usepackage[OT2,T1]{fontenc}
\usepackage{color}

\usepackage[
        colorlinks, citecolor=darkgreen,
        backref,
        pdfauthor={Elvira Lupoian},
]{hyperref}

\usepackage{comment}
\usepackage{multirow}
\usepackage[table]{xcolor}

\usepackage{hyperref}
\hypersetup{
    colorlinks=true,
    linkcolor=darkgreen,
    filecolor=magenta,      
    urlcolor=darkgreen,
    pdftitle={Overleaf Example},
    pdfpagemode=FullScreen,
    }
    
\newcommand{\Q}{\mathbb{Q}}
\newcommand{\R}{\mathbb{R}}
\newcommand{\Z}{\mathbb{Z}}
\newcommand{\Qbar}{{\overline{\Q}}}

\DeclareMathOperator{\Gal}{Gal}

\newcommand{\isom}{\simeq}

\numberwithin{equation}{section}

\newtheorem{thm}{Theorem}

\theoremstyle{definition}
\newtheorem*{defn}{Definition}

\newtheorem*{Ack}{Acknowledgements}

\theoremstyle{remark}	  
\newtheorem{rem}{Remark}

\newtheorem*{c1}{Case 1}
\newtheorem*{c2}{Case 2}

\theoremstyle{definition}

\theoremstyle{remark}

\definecolor{darkgreen}{rgb}{0,0.5,0}

\DeclareRobustCommand{\SkipTocEntry}[5]{}

\begin{document}

\title{Two-Torsion Subgroups of Some Modular Jacobians}

\begin{abstract}
We give a practical method to compute the 2-torsion subgroup of the Jacobian of a non-hyperelliptic curve of genus $3$, $4$ or $5$. The method is based on the correspondence between the 2-torsion subgroup and the theta hyperplanes to the curve. The correspondence is used to explicitly write down a zero-dimensional scheme whose points correspond to elements of the $2$-torsion subgroup. Using $p$-adic or complex approximations (obtained via Hensel lifting or homotopy continuation and Newton-Raphson)  and lattice reduction we are then able to determine the points of our zero-dimensional scheme and hence the $2$-torsion points. We demonstrate the practicality of our method by computing the $2$-torsion of the modular Jacobians $J_{0}\left( N \right)$ for $N = 42, 55, 63, 72, 75$. As a result of this we are able to verify the generalised Ogg conjecture for these values.  
\end{abstract}

\author{Elvira Lupoian}

\address{Mathematics Institute\\
    University of Warwick\\
    CV4 7AL \\
    United Kingdom}

\email{e.lupoian@warwick.ac.uk}
\date{\today}
\thanks{The author is supported by the EPSRC studentship}
\keywords{Two Torsion, Modular Jacobians,  Generalised Ogg Conjecture}
\subjclass[2020]{11G30}
 
\maketitle

\section{Introduction}
Let $X$ be a smooth, projective genus $g \ge 1$ curve over $\Q$ and let $J$ be its Jacobian variety. The Mordell-Weil theorem states that the set of  $K$-rational points of $J$ is a finitely generated group for any number field $K$; that is, $J \left( K \right) \isom J\left( K \right)_{\text{tors}} \oplus \Z^{r}$ for some integer $r \geq 0$ and a finite group $J\left( K \right)_{\text{tors}}$. In this paper we will be concerned with the finite torsion subgroup $J\left( K \right)_{\text{tors}}$, more specifically with the $2-$torsion part $J \left( K \right) \left[ 2 \right ] = \{ D \in J \left( K \right) \ \vert \ 2D = 0 \}$. If $X$ is hyperelliptic, then it is easy to compute the $2$-torsion of its Jacobian, for example see \cite{BS2tors} or \cite{scha2tors}. 
We describe a method for calculating the entire $2-$torsion subgroup $J\left( \Qbar \right) \left[ 2 \right] \isom \left( \Z / 2\Z \right)^{2g}$ for $X$ non-hyperelliptic and $g = 3, 4 \ \text{and} \ 5$.

The theoretical basis of  this method is the well-known description of 2-torsion points on the Jacobian as the difference of two odd theta characteristics \cite{caporaso}. A theta characteristic to the curve is a linear equivalence class of a degree $g-1$ divisor on $X$ which when doubled is equal to the canonical class. The parity of a theta characteristic is simply the parity of the dimension of its Riemann-Roch space. As our curve is non-hyperelliptic, it is canonically embedded in projective space and therefore all odd theta characteristics are in bijective correspondence with the hyperplanes to $X$ which intersect the curve at $g-1$ (not necessarily distinct) points, each with multiplicity $2$; we will call such planes theta hyperplanes. These are important geometric invariants of curves and have been studied extensively. For instance, Caporaso and Sernesi showed in \cite{CapEdo2} that a plane quartic, $g=3$, is completely determined by its theta hyperplanes, a result which they then generalised to all general canonical curves of genus $g \geq 4$ in \cite{CapEdo}. Lehavi showed that a curve can be effectively reconstructed from its theta hyerplanes in the genus $3$, $4$ and $5$ cases (see \cite{Leh1}, \cite{Leh2} and \cite{Leh3} respectively). The method described in this paper can be used to compute the theta hyperplanes to $X$ and thus the $2$-torsion points of its Jacobian, which are simply linear equivalence classes of differences of the divisors formed by intersecting the theta hyperplanes with the curve and multiplying by $\frac{1}{2}$. In Section $3$, we will construct zero-dimensional schemes whose points correspond to the theta hyperplanes to $X$. The points of such a zero-dimensional scheme are usually defined over a number field of fairly large degree, especially in the genus $4$ and $5$ cases, and as a result are usually impractical to compute using Gr\"{o}bner bases. In Sections $4$ and $5$ we describe how such points can be computed, firstly by approximating and then by searching for short vectors in an appropriate lattice. The points are approximated either $p$-adically, first by searching for non-singular points over a finite field and then lifting these points using Hensel's lemma, or as complex points, using homotopy continuation to obtain initial complex approximations and then applying Newton-Raphson to obtain complex approximations accurate to many decimal places. In Section 5 we use the approximations to define lattices in which we search for short vectors in order to find precise expressions for the points of our schemes. Using lattice reduction to find algebraic dependence is a standard technique, for instance it is described in \cite[Section 2.7.2]{cohen}. In Section $6$ we compute the $2$-torsion subgroup in some explicit examples. 

The main motivation for this work was to verify the generalised Ogg conjecture for some values of $N$, previously known up to $2$-torsion. Let $N$ be a positive integer and denote by $J_{0}\left( N \right)$ be the Jacobian variety of the modular curve $X_{0}\left( N \right)$. Denote by $C_{0}\left( N \right)$ the cuspidal subgroup of $J_{0}\left( N \right) \left( \overline{\mathbb{Q}} \right)$; that is, the subgroup generated by classes of differences of cusps, and write $C_{0}\left( N \right)\left( \mathbb{Q} \right)$ for the subgroup of $C_{0}\left( N \right)$ stable under the action of $\Gal \left( \mathbb{\overline{Q}} / \mathbb{Q} \right)$.  A consequence of the Manin-Drinfeld theorem \cite{dr}  is that $C_{0}\left( N \right) \left( \mathbb{Q} \right) \subseteq  J_{0}\left( N \right) \left( \mathbb{Q} \right)_{\text{tors}}$.  The generalised Ogg conjecture \cite{oggc} states that this is in fact an equality.  In \cite{OS} Ozman and Siksek proved this  for $N = 34, 38, 44, 45, 51, 52, 54, 56, 64, 81$. We compute the 2-torsion subgroup of $J_{0}\left( N \right)$ and using previous calculations of  \cite{OS},  we verify the generalised Ogg conjecture for $N = 42, 55, 63, 72, 75$, where $g=5$ for all these values of $N$.

\begin{thm}
The generalised Ogg conjecture holds for $N = 42, 55, 63, 72, 75$.
\end{thm}
The explicit generators of the $2$-torsion subgroups of the above, as well as the \texttt{Magma} code used to compute them can be found at 
\begin{center}
   \href{ https://github.com/ElviraLupoian/TwoTorsionSubgroups}{https://github.com/ElviraLupoian/TwoTorsionSubgroups}
\end{center}
Algorithms to compute theta hyperplanes to canonical curves have also been developed by 
Bruin, Poonen and Stoll \cite[Section 12]{BPS} in the case of $g = 3$;  where the $2$-torsion subgroup is used for explicit computations of  Selmer groups, and an algorithm to compute the tritangents, in the genus $4$ case is described by Stoll in \cite[Section 4]{StollTriTan}.

Methods to compute torsion points on Jacobians using Hensel lifting are also described by Dokchitser and Doris \cite{dokchitser}, and Mascot \cite{mascot}. We also use Hensel lifting in the genus $3$ and $4$ examples.  For our genus 5 examples, we found it impractical to search for nonsingular points over finite fields that can be used as initial approximations in our Hensel lifting. Instead, we use homotopy continuation \cite{ver}, to obtain initial complex approximations to the points of our zero-dimensional scheme. These are then lifted to very high precision, approximately $2000$ decimal places, using Newton-Raphson. These high precision approximations were crucial for finding the algebraic expressions of the theta hyperplanes. 

\begin{Ack}I would like to thank my supervisors Samir Siksek and Damiano Testa for their continued support and the many valuable conversations throughout this project, and Filip Najman for finding the mistake in the original proof of Theorem 1. I would also like to thank the anonymous referee for their detailed feedback, their suggestions and corrections have improved this paper greatly. 
\end{Ack}

\section{Preliminaries}
\label{Pre}
Let $X$ be a smooth, projective complex curve of genus $g$.
\begin{defn}
A $\mathbf{theta \ characteristic}$ on $X$ is a degree $g-1$ divisor class $D \in \text{Pic}\left( X \right)$ such that $2D = K_{X}$,  where $K_{X}$ is the canonical class on $X$.
The $\mathbf{parity}$ of a theta characteristic $D$ is the parity of $h^{0}\left( X, D \right)$.
\end{defn}

Let $J$ be the Jacobian variety of $X$. For any two theta characteristics $D_{1}, D_{2}$, the equivalence class of the difference $D_{1} - D_{2}$ is a $2$-torsion point of $J$, so there are $2^{2g}$ theta characteristics. It is well known that there are precisely $2^{g-1}\left( 2^{g} - 1\right)$ odd theta characteristics and $2^{g-1}\left( 2^{g} + 1  \right)$ even theta characteristics (see \cite[Chapter 5]{dog}).

The following result, stated in \cite{caporaso}, is a clear consequence of results of \cite[Chapter 5]{dog} and it is also proved in \cite[Cor 5.3]{BPS} for genus $g \geq 2$.
\begin{thm}
The $2$-torsion subgroup of $J$ is generated by differences of odd theta characteristics. 
\end{thm}

As described in \cite{caporaso}, theta characteristics can also be interpreted geometrically, as we now explain.
\begin{defn}
A $\mathbf{theta \ hyperplane}$ to $X$ is a hyperplane $H \subset \mathbb{P}^{g-1}$ tangent to $X$ at $g-1$ points.
\end{defn}

Suppose further that $X$ is non-hyperelliptic and it is embedded in $\mathbb{P}^{g-1}$ by its canonical embedding. Let $H$ be a theta hyperplane to $X$. The intersection of this hyperplane and the curve gives the following divisor on $X$:
\begin{center}
$H \cdot X = \displaystyle \sum_{i=1}^{g-1} 2P_{i}$ 
\end{center} 
for some $P_{i} \in X$. The equivalence class of the divisor $\frac{1}{2} H \cdot X$ is an effective theta characteristic. Conversely, given an effective theta characteristic, there exists a hyperplane $H \subset \mathbb{P}^{g-1}$ such that $D = \frac{1}{2}H \cdot X$, so $H$ is tangent to $X$ at $g-1 $ points. 
 
 For $X$ as above, odd theta characteristics coincide with the effective theta characteristics, and thus we can state the following result (see \cite{CapEdo}). 
\begin{thm}
There are  $2^{g-1}\left( 2^{g} -1 \right)$ theta hyperplanes to $X$, in natural bijection with the odd theta characteristics of $X$. 
 \end{thm}

 The $2-$torsion subgroup is computed as follows. Let $TH$ be the set of defining equations of the theta hyperplanes to $X$. We observe that the quotient of two elements of $TH$ is an element of the function field of $X$ and thus we can consider its divisor. The $2$-torsion subgroup of $J$ is generated by classes of divisors of the form 
 \begin{center}
 $\frac{1}{2} \left( \text{div} \left( \frac{a}{b} \right) \right)= \displaystyle \sum_{i=1}^{g-1} P_{i}  -     \displaystyle \sum_{j=1}^{g-1} Q_{j}$
 \end{center}
 for some $a,b \in TH$, where the $P_{i}$'s and $Q_{j}$'s are (not necessarily distinct) points on $X$. 

\section{Schemes of theta hyperplanes}
\label{Sch}
Let $X$ be a complete, nonsingular and non-hyperelliptic curve over an algebraically closed field $k$, of genus $g$.  The image of the canonical embedding of  $X$ into $\mathbb{P}^{g-1}$ is a curve of degree $2g - 2$ and this gives a model of the curve.

\begin{itemize}
\item $g = 3$ : the canonical embedding is a plane quartic in $\mathbb{P}^{2}$;
\item $g = 4$ : the canonical model is the intersection of a quadric and a cubic surface in $\mathbb{P}^{3}$;
\item $g = 5$ : the canonical model is the intersection of 3 quadrics in $\mathbb{P}^{4}$.
\end{itemize}
From now on, we will assume that $X$ has a model over $\mathbb{Q}$.  
In this section we construct a zero-dimensional scheme whose points correspond to theta hyperplanes to $X$, when $X$ has genus $3$, $4$ and $5$.

\subsection{The Genus 3 Case: Scheme of Bitangents}
Suppose $X$ has genus 3 and so the curve has a model:
\begin{center}
$X : f\left( x_{1}, x_{2}, x_{3} \right) = 0 $, 
\end{center}
where $f \in \mathbb{Q} \left[ x_{1}, x_{2}, x_{3} \right]$ is a homogeneous polynomial of degree 4.  A theta hyperplane to $X$ is a plane intersecting the curve in 2 points and we call such hyperplanes bitangent lines. 

We work on an affine chart and de-homogenise with respect to the appropriate coordinates.  For notation purposes, we work on the affine chart $\lbrace x_{3} = 1 \rbrace $. A bitangent to the affine curve is given by a polynomial 
\begin{center}
$b_{1}x_{1} + b_{2}x_{2} + b_{3} = 0 $,
\end{center}
for some $b_{i} \in \mathbb{\overline{Q}}$ with $ \left( b_{1},b_{2} \right) \neq \left( 0 , 0 \right)$. Suppose $b_{2} \ne 0 $, so rescaling and rearranging gives 
\begin{center}
$x_{2} = a_{1}x_{1} + a_{2} $,
\end{center}
for some $a_{1}, a_{2} \in \mathbb{\overline{Q}}$. The intersection of the affine curve with this line is described by 
\begin{center}
$F(x_{1} ) = f\left( x_{1} , a_{1} x_{1} + a_{2},  1 \right),$
\end{center}
where $F \in \mathbb{Q} \left[ a_{1}, a_{2}\right] \left[ x_{1} \right] $ has degree 4,  and $F$ is necessarily a square if the given line is a bitangent.  Equivalently, there exist $a_{3}, a_{4} \in \mathbb{\overline{Q}} $ such that 
\begin{center}
$F \left( x_{1} \right) = l \left( x_{1}^{2} + a_{3}x_{1} + a_{4} \right)^{2}$
\end{center}
where $l$ is the coefficient of $x_{1}^{4}$ in $F$. Equating coefficients in the above expression gives $4$ equations $e_{1}, \ldots , e_{4}$ in $a_{1}, \ldots, a_{4}$ which define a zero-dimensional scheme $S$,  whose points correspond to bitangents to the curve. 
\begin{rem}
The total number of bitangents to a plane curve is $28$ and hence the degree of $S$ is at most 28. If the degree is strictly less than 28, we can repeat the above,  working on a different affine chart or with lines of a different form to obtain all 28 bitangents to $X$.
\end{rem}

\subsection{The Genus 4 Case: Scheme of Tritangents}
Suppose $X$ has genus 4, and so a canonical model of the curve is the intersection of a quadric and a cubic.
\begin{center}
$X : f\left( x_{1}, x_{2}, x_{3}, x_{4} \right) = g\left( x_{1}, x_{2}, x_{3}, x_{4} \right) =0 $
\end{center}
where $f, g \in  \mathbb{Q} \left[ x_{1}, x_{2}, x_{3}, x_{4} \right]$ are homogeneous of degree $2$, $3$ respectively.  A theta hyperplane to $X$ is a plane intersecting the curve in $3$ (not necessarily distinct) points, each with multiplicity two and we will call such hyperplanes tritangents planes.

We work on an affine chart, say $ \lbrace x_{4} = 1 \rbrace $ and similar to above,  we can assume that some tritangent planes to the affine curve are cut out by equations of the form
\begin{center}
$x_{3} = a_{1}x_{1} + a_{2}x_{2} + a_{3} $,
\end{center}
for some $a_{1}, a_{2}, a_{3} \in \mathbb{\overline{Q}}$. The intersection of the affine curve with the plane is described by the two expressions $F$ and $G$,
\begin{align*}
F\left( x_{1}, x_{2} \right) & \  =  \ f \left( x_{1},  x_{2},  a_{1}x_{1} + a_{2}x_{2} + a_{3} , 1 \right), \\
G \left( x_{1}, x_{2} \right) &  \  =  \ g \left( x_{1}, x_{2},  a_{1}x_{1} + a_{2}x_{2} + a_{3} , 1 \right),
\end{align*}
which can be viewed as polynomials in $x_{2}$ with coefficients in $\mathbb{Q} \left[ a_{1}, a_{2}, a_{3}, x_{1} \right]$. Taking the resultant of $F$ and $G$ gives a polynomial $R$ in $x_{1}$ with coefficients in $\mathbb{Q} \left [ a_{1}, a_{2}, a_{3} \right]$,
\begin{center}
$R \left( x_{1} \right) = \text{Res}\left( F, G, x_{2} \right) \in \mathbb{Q} \left[ a_{1}, a_{2}, a_{3} \right]\left[ x_{1} \right]  $.
\end{center}
This is a degree 6 polynomial; and it is necessarily a square if the given plane is a tritangent. That is, there exist $a_{4}, a_{5}, a_{6} \in \mathbb{\overline{Q}}$ such that 
\begin{center}
$R \left( x_{1} \right) = l \left( x_{1}^{3} + a_{4}x_{1}^{2} + a_{5}x_{1} + a_{6} \right)^{2},$
\end{center}
where $l$ is the coefficient of $x_{1}^{6}$ in $R$. The coefficients of the above expression $e_{1}, \ldots, e_{6}$ are polynomials in $a_{1}, \ldots, a_{6}$ with rational coefficients. To define the scheme of tritangents we require the additional equation 
\begin{center}
$e_{7} : a_{7} \cdot \Delta + 1 = 0$,
\end{center}
where $\Delta$ is the discriminant of $x_{1}^{3} + a_{4}x_{1}^{2} + a_{5}x_{1} + a_{6}$. This equation ensures that $\Delta \neq 0 $ and avoids singularities on the scheme. The seven equations derived above define a zero-dimensional scheme $S$, whose points correspond to tritangent planes to $X$.

\begin{rem}
The total number of tritangent planes to such a curve is $120$ and hence the degree of $S$ is at most $120$. As in the genus $3$ case, if the degree is strictly less than $120$, we can repeat the above, working on a different affine chart or with planes of a different form to obtain all $120$ tritangent planes to $X$.
\end{rem}

\subsection{The Genus 5 Case: Scheme of Quadritangents}
\label{Sch3}
Suppose $X$ has genus $5$ and so the curve has a canonical model of the following form:
\begin{center}
$ X : f_{1}\left( x_{1},  \ldots ,   x_{5} \right) = f_{2}\left( x_{1}, \ldots, x_{5} \right) = f_{3}\left( x_{1}, \ldots, x_{5} \right) = 0,$
\end{center}
where $f_{i} \in \mathbb{Q} \left[ x_{1}, x_{2}, x_{3}, x_{4}, x_{5} \right]$ are homogeneous of degree $2$. A theta hyperplane to $X$ is a hyperplane intersecting the curve in $4$ (not necessarily distinct) points, each with multiplicity $2$; and such hyperplanes will be called quadritangent planes throughout this paper.
 
 As in the previous two cases, to construct the scheme of quadritangent planes we intersect the affine curve with a  given plane and eliminate variables, one at a time, until we obtain an expression in one variable which is necessarily a square if the given plane is a quadritangent.  This is dependent on the model of the curve.  In this subsection we give a general idea of how such a scheme can be defined; for an explicit example see Section 6.3.
With $X$ as above,  we work on an affine chart,  say $\lbrace x_{5} = 1 \rbrace$, and as in previous cases we can assume that some quadritangents to the curve are given by equations of the form:
 \begin{center}
 $x_{4}  = a_{1}x_{1} + a_{2}x_{2} + a_{3}x_{3} + a_{4}$
 \end{center}
 for some $a_{1}, \ldots, a_{4} \in \Qbar$. The intersection of the affine curve with the above plane is described by $3$ polynomials:  
 \begin{center}
 $F_{i}\left( x_{1}, x_{2}, x_{3} \right) = f_{i}\left( x_{1}, x_{2}, x_{3},  a_{1}x_{1} + a_{2}x_{2} + a_{3}x_{3} + a_{4}, 1 \right)$,  \ \   $i = 1, 2, 3$,
 \end{center}
 in $x_{1}, x_{2}, x_{3}$ with coefficients in $\mathbb{Q} \left[ a_{1}, a_{2}, a_{3}, a_{4} \right]$ and monomials of total degree at most $2$.
 We now eliminate one of the variables.
\begin{c1} Suppose one of the $F_{i}$, say $F_{1}$, is linear in one of the variables, say $x_{3}$, so $F_{1}$ can be written as 
\begin{center}
$F_{1}\left( x_{1}, x_{2}, x_{3} \right) \ = \  g_{1}\left(x_{1}, x_{2} \right) x_{3} + g_{2}\left( x_{1}, x_{2} \right),$
\end{center}
where $g_{1}, g_{2} \in \mathbb{Q}\left[ a_{1}, a_{2}, a_{3}, a_{4} \right] \left[ x_{1}, x_{2} \right]$. If $g_{1}\left( x_{1}, x_{2} \right) \neq 0$, we obtain in expression for $x_{3}$,
\begin{center}
$x_{3} = - \frac{g_{2}\left( x_{1}, x_{2} \right)}{g_{1}\left( x_{1}, x_{2} \right)} = \varphi \left( x_{1}, x_{2} \right),$
\end{center}
which can be substituted into $F_{2}$ and $F_{3}$.
\end{c1}
\begin{c2}
If none of the $F_{i}$ are linear in one variable, then we can write $F_{1}$ and $F_{2}$ as quadratics in $x_{3}$:
\begin{align*}
F_{1} & = g_{1,2} x_{3}^{2} + g_{1,1} x_{3} + g_{1,0}, \\
F_{2} & = g_{2,2} x_{3}^{2} + g_{2,1}x_{3} + g_{2,0}, 
\end{align*}
where $g_{i,j} \in \mathbb{Q} \left[ a_{1}, a_{2}, a_{3}, a_{4} \right] \left[ x_{1}, x_{2} \right]$. If $g_{1,2} \ne 0$ and $g_{2,2} \ne 0 $ we cross multiply to obtain the following expression
\begin{align*}
g_{1,2}F_{2} - g_{2,2}F_{1} & = \left( g_{1,2}g_{2,1} - g_{2,2} g_{1,1} \right) x_{3}  + g_{1,2} g_{2,0} - g_{2,2} g_{1,0} = 0  \\
\end{align*}
and we use this to define $\alpha = g_{1,2}g_{2,1} - g_{2,2} g_{1,1}$ and $\beta =  g_{1,2} g_{2,0} - g_{2,2} g_{1,0} $. Additionally, if we assume that $\alpha \ne 0 $:
\begin{center}
$x_{3} = - \frac{\beta}{\alpha} = \varphi \left( x_{1}, x_{2} \right)$
\end{center}
and we obtain an expression for $x_{3}$ which can be substituted into the $F_{i}$'s.
\end{c2}

In deriving an expression for $x_{3}$, we assumed additional conditions, such as $g_{i,j} \neq 0 $.  These conditions can be converted into equations, which will be added to the list of equations when defining the scheme. For example, if we assumed $a_{1} \neq 1$, the corresponding equation will be $b_{1}\left( a_{1} - 1 \right) + 1 =0$ for a newly introduced variable $b_{1}$.

Given such an expression $x_{3} = \varphi \left( x_{1}, x_{2} \right)$, we can substitute this into the $F_{i}$'s to obtain:
\begin{center}
$H_{i} \left( x_{1} x_{2} \right) = F_{i} \left( x_{1}, x_{2}, \varphi \left( x_{1}, x_{2} \right) \right) $ \ $ i = 1, 2,3$,
\end{center}
where $2$ of these equations are independent, and clearing denominators we obtain polynomials $h_{1}, h_{2} \in \mathbb{Q} \left[ a_{1}, a_{2}, a_{3}, a_{4} \right] \left[ x_{1}, x_{2} \right]$, whose monomials have total degree at most 4.  We now eliminate one of the two remaining variables. One practical way in which this can done is by considering whether or not $h_{1}$ and $h_{2}$ are quadratics in one variable.
\begin{c1}
If $h_{1}, h_{2}$ are both quadratic in one variable, say $x_{2}$, and thus they can be written as:
\begin{align*}
h_{1}\left( x_{1}, x_{2} \right) & = h_{1,2}x_{2}^{2} + h_{1,1}x_{2} + h_{1,0}, \\
h_{2}\left( x_{1}, x_{2} \right) & = h_{2,2}x_{2}^{2} + h_{2,1}x_{2} + h_{2,0} 
\end{align*}
where $h_{i,j} \in \mathbb{Q} \left[ a_{1}, a_{2}, a_{3}, a_{4} \right] \left[ x_{1} \right]$. If $h_{1,2} \neq 0 $ and $h_{2,2} \neq 0 $, we cross multiply:
\begin{align*}
h_{2,2}h_{1} - h_{1,2}h_{2} = \left( h_{2,2} h_{1,1} - h_{1,2}h_{2,1} \right) x_{2} + h_{2,2}h_{1,0} - h_{1,2} h_{2,0} = 0.
\end{align*}
Define $\gamma = h_{2,2} h_{1,1} - h_{1,2}h_{2,1}$ and $\delta = h_{2,2}h_{1,0} - h_{1,2} h_{2,0} $. If $\gamma \neq 0 $ then: 
\begin{center}
$x_{2} = - \frac{-\delta}{\gamma} = \psi \left( x_{1} \right) $.
\end{center}
We substitute this into $h_{1}$ 
\begin{center}
$R\left( x_{1} \right) = H_{1} \left( x_{1}, \psi \left(x_{1} \right) \right)$
\end{center}
and clearing denominators, we obtain an expression $r \left( x_{1} \right) \in \mathbb{Q} \left[ a_{1}, a_{2}, a_{3}, a_{4} \right] \left[ x_{1} \right]$. One of the factors of $r$ is a degree 8 polynomial $h \left( x_{1} \right) \in  \mathbb{Q} \left[ a_{1}, a_{2}, a_{3}, a_{4} \right] \left[ x_{1} \right]$.
\end{c1} 

\begin{c2}
If $h_{1}, h_{2}$ are not both quadratic in the same variable, we take the resultant with respect to $x_{2}$:
\begin{center}
$R\left(x_{1} \right)  = \text{Res}\left( h_{1}, h_{2}, x_{2} \right) \in \mathbb{Q}\left[ a_{1}, a_{2}, a_{3},a_{4} \right] \left[ x_{1} \right]$
\end{center}
and as in the previous case, $R$ has a degree 8 factor $h\left( x _{1} \right) \in  \mathbb{Q} \left[ a_{1}, a_{2}, a_{3}, a_{4} \right] \left[ x_{1} \right] $.
\end{c2}
Taking the resultant can cut out a larger zero-dimensional scheme than the scheme of quadritangents. With notation as above, requiring $h$ to be a square will ensure that $x_{1}$ occurs with multiplicity 2, but the corresponding $x_{2}$ might not. In practice, we find that this is not a problem, as the additional points not corresponding to quadritangents can simply be discarded at the end.

If there are multiple possible polynomials $h$ resulting from the above, we simply repeat the procedure below for all possibilities. In our computations, we found that this was not an issue, in all examples there was only one choice of $h$.

If the given plane is a quadritangent then $h$ is necessarily a square. That is, there exist $a_{5}, a_{6}, a_{7},a_{8} \in \mathbb{\overline{Q}}$ such that:
\begin{center}
$h\left(x_{1} \right) = l \left( x_{1}^{4} + a_{5}x_{1}^{3} + a_{6}x_{1}^{2} + a_{7}x_{1} + a_{8} \right)^{2} $
\end{center}
where $l$ is the leading coefficient of $h$. The coefficients in the above expression give 8 equations $e_{1}, \ldots, e_{8}$ with rational coefficients in $a_{1}, \ldots, a_{8}$. To define a scheme of quadritangents we also require the additional equation 
\begin{center}
$e_{9} : a_{9} \cdot \Delta + 1 =0 $
\end{center}
where $\Delta$ is the discriminant of $x_{1}^{4} + a_{5}x_{1}^{3} + a_{6}x_{1}^{2} + a_{7}x_{1} + a_{8}$. This ensures that the discriminant is non-zero and avoids singularities on our scheme.  The above equations along with any equations arising from any conditions necessary to derive $h$ define a zero-dimensional scheme whose points correspond to quadritangent planes to the curve. 

\section{Approximate Theta Hyperplanes }
\label{apr}
In the previous section we described a method of deriving equations which define zero-dimension schemes whose points correspond to coefficients of the theta hyperplanes to the curve. Gr\"{o}bner basis techniques can in theory be used to compute the points of such a zero-dimensional scheme. For example, the following two \texttt{Magma} commands do precisely this: \texttt{PointsOverSplittingField} and \texttt{Points}. The input for the former is a set of equations defining a zero-dimensional scheme and its output is the solution set of the system of equations. The latter command is less ambitious. It is designed to give the set of $K$-rational points of a zero- dimension scheme $S$, where $K$ is the field of definition of $S$. We found that \texttt{PointsOverSplittingField} is extremely slow in our examples, and in fact we were not able to use this even in the genus $3$ case.  Given the field of definition of the theta hyperplanes, the command \texttt{Points} was sometimes successful in determining the points, most notably in the genus $4$ example presented in Section 6, were it took around 10 minutes to compute the points given the degree $36$ number field over which all tritangents to the curve are defined. The field of definition of the tritangent planes can be determined using the method described in the two sections which follow. However, this command is still inefficient in the genus $5$ case, and we needed to use the methods of the Sections $4$ and $5$.

In this section we describe two methods of approximating points on a scheme of theta hyperplanes $S$ to a curve $X$.

\subsection{$p$-adic Approximations}
We briefly describe a method of approximating the points of our scheme $p$-adically. This method is also used in \cite{dokchitser} and \cite{mascot}.

Suppose $p$ is an odd prime of good reduction for the curve. We view the defining equations of $S$ over the finite field $\mathbb{F}_{p}$ and search for points of $S$ over an extension $\mathbb{F}_{p^{n}}$ for some $n \in \mathbb{N}$ . Any smooth points in $S\left( \mathbb{F}_{p^{n}} \right)$ are lifted using the following multivariate version of Hensel's Lemma (see \cite{conradmulthensel}).

\begin{thm}{(Multivariate Hensel Lemma) }
Let $K$ be a non-archimedean local field, $\nu_{K}$ a valuation on $K$ and $\mathcal{O}_{K}$ its ring of integers.  Suppose $F = \left( F_{1}, \ldots, F_{m} \right) \in \mathcal{O}_{K}\left[ x_{1}, \ldots, x_{m} \right]$ and $P = \left( P_{1}, \ldots, P_{m} \right) \in \mathcal{O}_{K}^{m}$; and let $v_{1} = \text{min} \lbrace  \nu_{K}\left( P_{i} \right) : i  = 1, \ldots, m \rbrace $ and $v_{2} = \nu_{K} \left( D_{F}\left(P \right) \right)$, where $D_{F}$ denotes the determinant of the Jacobian matrix of $F$. If $v_{1} > 2v_{2}$, then there exists a unique $Q = \left( Q_{1}, \ldots, Q_{m} \right) \in \mathcal{O}_{K}^{m}$ such that $F \left( Q \right) = 0 $ and $\nu_{K}\left( P_{i} - Q_{i} \right) \geqslant v_{1} - v_{2} $ for all $i$.
\end{thm}

 The proof of Hensel's lemma is constructive in the sense that for a  point $x \in S $ which satisfies the hypothesis of the above theorem, we can  find a prime ideal $\mathfrak{p}$ in  $\mathcal{O}_{K}$,  where $K$ is some number field and $\mathcal{O}_{K}$ is its ring of integers, and construct  a sequence $ \lbrace a_{i} \rbrace_{i \geq 1 } \subset \mathcal{O}_{K}$ satisfying the following:
\begin{itemize}
\item $a_{k} \equiv x \ \text{mod} \  \mathfrak{p}^{k}$ for all $k \geq 1$,
\item $a_{k} \equiv a_{k-1} \ \text{mod} \ \mathfrak{p}^{k-1}$ for all $k \geq 2$,
\item $E_{S}\left( a_{k} \right) \equiv  \ 0  \ \text{mod} \  \mathfrak{p}^{k}$ for all $k \geq 1$,
\end{itemize}
where $E_{S} = \left( e_{1}, \ldots, e_{n} \right)$ are the defining equations of $S$. Note that the congruence in the last expressions is understood to be coordinate wise. 

\begin{rem}
For large degree schemes, large $p$ and $n$ are required to find sufficiently many smooth points. In particular, this method was extremely inefficient in the computations required for the proof of Theorem 1; complex approximations were used for those calculations. 
\end{rem}

\subsection{Complex Approximations}
\label{apr2}
The points of $S$ can also be approximated as complex points using the Newton-Raphson method.  We give a brief overview of this,  a detailed explanation can be found in \cite[Chapter 5]{NA}.  

Let $E = \left( e_{1}, \ldots, e_{n} \right)$ be the defining equations of $S$.  We view $E$ as a function $ \mathbb{C}^{n} \longrightarrow \mathbb{C}^{n}$. Let $dE$ be the Jacobian matrix of $E$ and suppose  $\mathbf{x}_{0}$ is an approximate solution to $E$ with $dE\left( \mathbf{x}_{0} \right)$ invertible.  For $k \geq 1$, define 
\begin{center}
$\mathbf{x}_{k} = \mathbf{x}_{k-1} - dE\left( \mathbf{x}_{k-1} \right)^{-1}E\left( \mathbf{x}_{k-1} \right)$
\end{center}
Provided the initial approximation $\mathbf{x_{0}}$ is a good approximation,  the resulting sequence $\lbrace \mathbf{x}_{k} \rbrace_{k \geq 0 }$ converges to a solution of $E$, with each iterate having increased precision. In fact, at each step the number of decimal places to which the approximation is accurate roughly doubles \cite[Section 5.8]{NA}. This was the method used to compute approximations of the theta hyperplanes of the curves stated in Theorem 1.  

This method requires good initial complex approximations to the solutions of $E$. These can be obtained using Homotopy Continuation and its implementation in the numerical analysis package \texttt{Julia}.

\subsubsection{Homotopy Continuation}
Homotopy continuation is a method for numerically approximating the solutions of a system of polynomial equations by deforming the solutions of a similar system whose solutions are known.  We give a brief sketch of the idea,  but a more detailed explanation of this theory can be found in \cite{ver} or \cite{NA}.
                
The total degree of $E$ is defined as $\text{deg} \left( E \right) = \displaystyle \prod_{i=1}^{n} \text{deg} \left( e_{i} \right)$, where $\text{deg} \left( e_{i} \right)$ is the maximum of the total degrees of the monomials of $e_{i}$.

Let $F$ be a system of $n$ polynomials in $ \mathbf{a} = \left( a_{1}, \ldots, a_{n} \right)$,  which has exactly $\text{deg}\left( E \right)$ solutions and these solutions are known.  The system $F$ will be known as a start system.  The standard homotopy of $F$ and $E$ is the function 

\begin{align*}
 & H :  \mathbb{C}^{n} \times \left[ 0, 1 \right] \longrightarrow \mathbb{C}^{n} \\
    & H \left( \mathbf{a}, t \right) \  = \  \left( 1 - t \right) F\left( \mathbf{a} \right) + t E \left( \mathbf{a} \right)
\end{align*}
Fix $N \in \mathbb{N}$, and for any $s \in \left[ 0, N \right] \cap \mathbb{N}$ define $H_{s}\left( \mathbf{a} \right) = H \left( \mathbf{a}, s/N \right)$, this is a system of $n$ polynomials in $ \mathbf{a} = \left( a_{1}, \ldots, a_{n} \right)$.

For $N$ large enough, the solutions of $H_{s} \left( \mathbf{a} \right)$ are good initial approximations of the solutions of $H_{s+1}\left( \mathbf{a} \right)$, and using the Newton-Raphson method we can greatly increase their precision.  The solutions of $H_{0} \left( \mathbf{a} \right) = F\left( \mathbf{a} \right)$ are known, and they can be used to define solution paths to approximate solutions of $H_{N}\left( \mathbf{a} \right) = E \left( \mathbf{a} \right)$.

There are two important things to highlight.
\begin{itemize}
\item[1.] Given any $E$,  a start system (and its solutions) can always be computed.
\item[2.] A start system can be modified to ensure solutions paths do not cross and converge to approximate solutions of $E$.
\end{itemize}

Homotopy continuation is implemented in the \texttt{Julia} with the simple command \texttt{solve}, whose input is any complex system of equations which has a finite number of solutions (see \cite{JulHC}). 

\begin{rem}
The implementation of homotopy continuation in \texttt{Julia} gives approximates to solutions of $E$ which are accurate to 16 decimal places.  For the computations required for the proof of Theorem 1, we used the approximate solutions and 600 iterations of Newton-Raphson to obtain an accuracy of 2000 decimal places,  which was sufficient for those calculations. 
\end{rem}

\section{Algebraic Theta Hyperplanes}
\label{alg}
Given a scheme of theta hyperplanes $S$, as in Section $3$, and an approximation $\mathbf{a} = \left( a_{1}, \ldots, a_{n} \right)$, as in Section $4$, to a point $\mathbf{x} = \left( x_{1}, \ldots, x_{n} \right)$ of $S$, we use lattice reduction to compute the minimal polynomials of the $x_{i}$ and hence determine $\mathbf{x}$ exactly. The use of lattices and lattice reduction is a standard method when looking for minimal polynomials, for instance see \cite[Section 2.7.2]{cohen} where Cohen explains how the LLL-algorithm can can be used to determine linear and algebraic dependence. Similarly, Smart \cite[Chapter 6.1]{smart} shows that the LLL-algorithm can be used to find a polynomial which has a root that is a good approximation of $\pi$.

\subsection{Algebraic Hyperplanes and $p$-adic Approximations}
Suppose that $\mathbf{x} = \left( x_{1}, \ldots, x_{n} \right)$ is a point of $S$ for which we have a $p$-adic approximation $\mathbf{a} = \left( a_{1}, \ldots, a_{n} \right)$, and we wish to compute the minimal polynomials of the $x_{i}$. As in Section 4.1 there exist a number field $K$, a prime ideal $\mathfrak{p}$ of $\mathcal{O}_{K}$ and a point $ \left(  a_{1}, \ldots, a_{n} \right)$, with $a_{i} \in \mathcal{O}_{K}$ such that for some $k \in \mathbb{N}$,  $x_{i} \equiv a_{i} \ \text{mod} \ \mathfrak{p}^{k}$ for all $i \ge 1$.  We will call this $k$ the precision of this approximation; and as explained in Section 4.1, approximations with arbitrarily large precision $k$ can be computed.

The coefficients of the equations defining the theta hyperplanes are algebraic numbers and we can calculate their minimal polynomials by searching for short vectors in an appropriately defined lattice. Short vectors can be computed using the  Lenstra-Lenstra-Lov\'{a}sz lattice reduction algorithm, or LLL for short, which given a lattice, returns a reduced basis whose vectors have small norms, see \cite[Section 2.6]{cohen}. 

\subsubsection*{Minimal Polynomials}
Fix $i$, $1 \leq i \leq n$ and let $\theta = x_{i}$. In this subsection we describe a method to compute the minimal polynomial of $\theta$ using the $p$-adic approximations.  As $\theta$ is algebraic, there exist $d_{\theta} \in \mathbb{N}$ and $c_{0}, \ldots, c_{d_{\theta}} \in \mathbb{Z}$ such that 
\begin{center}
$c_{d_{\theta}} \theta^{d_{\theta}} + \cdots + c_{1} \theta + c_{0} = 0 $
\end{center}
 where $c_{d_{\theta}} \neq 0 $. As $\theta \equiv a_{i} \ \text{mod} \ \mathfrak{p}^{k}$,
 \begin{center}
 $c_{d_{\theta}} a_{i}^{d_{\theta}} + \cdots + c_{1} a_{i} + c_{0} \equiv 0 \ \text{mod} \ \mathfrak{p}^{k} $.
 \end{center}
Define the homomorphism 
\begin{center}
$\phi : \mathbb{Z}^{d_{\theta}+1} \longrightarrow \mathcal{O}_{K} / \mathfrak{p }^{k}$ \\
\medskip
$ \left( u_{0}, \ldots , u_{d_{\theta}} \right) \longmapsto  u_{d_{\theta}}a_{i}^{d_{\theta}} + \cdots + u_{1}a_{i} + u_{0} \ \text{mod} \ \mathfrak{p}^{k} $,
\end{center}
and let $\mathcal{L} = \text{ker} \left( \phi \right) \leqslant  \mathbb{Z}^{d_{\theta}+1}$. This is a discrete subgroup of $\mathbb{Z}^{d_{\theta}+1}$, which contains all elements of the form $P\mathbf{e}_{i}$, where $P$ is any integer which reduces to $0$ modulo $\mathfrak{p}^{k}$ and $\mathbf{e}_{i}$ is the $ith$ element in the standard orthonormal basis of $\mathbb{Z}^{d_{\theta}+1}$, and thus $\mathcal{L}$ is a full rank lattice in $\mathbb{Z}^{d_{\theta}+1}$.

Observe that the above lattice $\mathcal{L}$ can be constructed for any precision $k$, and $\left( c_{0}, \ldots , c_{d_{\theta}} \right) \in \mathcal{L}$ for any $k$. As $k$ increases, in general the average length of vectors in $\mathcal{L}$ increases and so eventually $\left( c_{0}, \ldots , c_{d_{\theta}} \right) $  should be the shortest vector in $\mathcal{L}$, or one of the shortest vectors. The length of a vector in a lattice is defined as the square root of its norm. To search for short vectors in $\mathcal{L}$ we use the \texttt{Magma} command \texttt{ShortestVectors} which given a lattice returns a sequence containing all vectors of the lattice which have the minimum non-zero norm, see \cite{MR1484478}. The \texttt{Magma} algorithm for computing short vectors first computes a reduced basis of the lattice using its efficient implementation of the LLL-algorithm, which is based on both Nguyen and Stehl\`{e}'s floating-point LLL algorithm \cite{NSLLL} and de Weger's exact integral algorithm \cite{dWLLL}; and then uses a closest vector algorithm, such as the one described by Fincke and Pohst \cite{FP}, to determine the shortest non-zero vector in $\mathcal{L}$. 

To determine what precision makes the shortest vector in $\mathcal{L}$ a suitable candidate for the coefficients of the minimal polynomial we use Hermite's theorem. 
\begin{thm}{(Hermite)}
There exist constants $\mu_{n} \in \mathbb{R}_{\ge 0}$ such that for any $n$-dimensional lattice $\mathcal{L}$ we have
\begin{center}
$M^{n} \leq \mu_{n} \Delta \left( \mathcal{L} \right)^{2} $
\end{center}
where $M$ is the length of the shortest non-zero vector in the lattice and $\Delta \left( \mathcal{L} \right)$ is the determinant of the lattice, as defined in \cite[Page 66]{smart}.
\end{thm}
There are bounds on these $\mu_{n}$ given in \cite[Page 66]{smart}. 
\begin{proof}
See \cite[ Page 66]{smart}
\end{proof}

In our case, we expect the determinant $\Delta \left( \mathcal{L} \right)$ to be around $I^{1/(d_{\theta} + 1)}$, where $I$ is the index $\left[ \mathbb{Z}^{d_{\theta} + 1 } : \mathcal{L} \right]$, and the homomorphism $\phi$ to be surjective. Thus, 
\begin{center}
    $I = \frac{\vert \left( \mathcal{O}_{K} / \mathfrak{p}^{k} \right)^{\times} \vert }{\vert \left( \Z / p^{k} \Z \right)^{\times} \vert } \approx \frac{\text{Norm}\left( \mathfrak{p} \right)^{k}}{p^{k}} $
\end{center}
and given our data, the last quantity can be explicitly computed. Therefore, if we are searching for a minimal polynomial whose coefficients are of approximate order $10^{n}$, then we require $k$ such that 
\begin{center}
   $ \left( \frac{\text{Norm} \left( \mathfrak{p}\right) }{p} \right)^{k/( d_{\theta} + 1 )} \ge \sqrt{\left( d_{\theta} + 1 \right)} \cdot 10^{n}$.
     \end{center}

For large  $k$, the length of the shortest non-zero vector in $\mathcal{L}$ should be significantly smaller than $\left(\Delta \left( \mathcal{L} \right)\right)^{\frac{1}{d_{\theta} + 1 }}$, and this would make it a suitable candidate for the vector whose entries are the coefficients of the minimal polynomial. 

In the construction above, we can replace $d_{\theta}$ by any positive integer $d$ and we search for short vectors in the lattice constructed with respect to this $d$. If our candidate $d$ is equal to or slightly bigger than $d_{\theta}$, then the following should hold.
\begin{itemize}
\item The length of $v_{\text{min}}$, the shortest non-zero vector in the lattice,  should  be significantly smaller than $ d \left(\mathcal{ L} \right)^{\frac{1}{d+1}}$ for a large enough $k$.
\item The polynomial $f_{\text{min}}$ whose vector of coefficient is $v_{\text{min}}$ should be irreducible or its factorisation should contain an irreducible polynomial of degree close to the degree of $f_{\text{min}}$.
\item The minimal polynomials should respect the Galois action on the theta hyperplanes, so multiple points should have the same minimal polynomial.
\end{itemize}

To summarise, the strategy for finding the coefficients of the minimal polynomial of $\theta$ is as follows.
\begin{itemize}
\item[1.] We start with a candidate for the degree $d$.
\item[2.] Define the homomorphism $\phi$ and the lattice $\mathcal{L}$, depending on degree $d$ and an approximation of precision $k$. Note that $k$ also depends on $d$.
\item[3.] In  $\mathcal{L}$ search for vectors which are shorter than, say $1/1000 \Delta\left( \mathcal{L} \right)^{\frac{1}{d+1}}$. 
\item[4.] If such a vector exists, verify that it satisfies the conditions stated above, and if they are all satisfied, it's extremely likely that this vector represents the coefficients of the minimal polynomial. Otherwise, return to 1 and replace $d$ by $d +1$.
\end{itemize}

The initial candidate for the degree is $d=1$, and we run through $d \in \mathbb{N}$, until a suitable vector is found. 
\subsubsection*{Coefficient Relations}
Once we have candidates for the minimal polynomials for all $ x_{1}, \ldots, x_{n} $ we want to identify which roots of the minimal polynomials correspond to the point approximated. This is especially important when our computations use hyperplanes whose defining equations have coefficients whose minimal polynomials have large degrees, since in these cases it's impractical to simply run through all tuples of possible roots and test which  define theta hyperplanes. 

 When the degrees are small, the simplest way of choosing the corresponding roots is through factorisation.  Suppose that we want to find the root of $f_{1}$ approximated by $a_{1}$. We  factor  $f_{i}$ over a number field $\tilde{K}$ over which $f_{1}$ splits into linear factors and for which $K$ is a subfield.  Let $\mathfrak{P}$ be a prime of $\tilde{K}$ over $\mathfrak{p}$. We factor $f_{i}$ over $\tilde{K}$,
\begin{center}
$f_{i} = s_{1}\ldots s_{d_{i}}$
\end{center}
where $s_{j}$ are linear factors. Then for some $k$, 
\begin{center}
$s_{k} \left( a_{1} \right) \equiv 0 \ \text{mod} \ \mathfrak{P}^{k}$
\end{center}
and thus $s_{k}$ corresponds to the root of $f_{1}$ approximated by $a_{1}$. 

When factorising is not efficient, for instance when the degrees of our polynomials are large, we may identify the required roots using lattice reduction. We search for relations amongst the coefficients to identify the required roots. Suppose that all our coefficients $x_{1},\ldots, x_{n}$ are contained in the number field $K = \Q \left( x_{1} \right)$ defined by $f_{1}$. Then, there exist $b_{0}, \ldots b_{d_{1}} \in \Z $ such that 
\begin{center}
$b_{d_{1}} x_{2} +  b_{d_{1} -1} x_{1}^{d_{1} - 1} + \ldots + b_{1} x_{1} + b_{0} =0 $
\end{center}
where $d_{1}$ is the degree of $f_{1}$.
Define the homomorphism,
\begin{center}
$r_{k} : \mathbb{Z}^{d_{1} + 1 } \longrightarrow \mathcal{O}_{K}/ \mathfrak{p}^{k}$\\
\medskip
$\left( u_{0}, \ldots, u_{d_{1} -1}, u_{d_{1}} \right) \longmapsto u_{0} + u_{1}a_{1} + \ldots + u_{d_{1} -1}a_{1}^{d_{1} -1} - u_{d_{1}}a_{2} \  \text{mod} \ \mathfrak{p}^{k}.$
\end{center}
Let $R_{k} = \text{ker}\left( r_{k} \right) < \mathbb{Z}^{d_{1} + 1 }$. This is a full rank lattice in $ \mathbb{Z}^{d_{1} + 1 }$ and arguing as before, for $k$ large enough $\left( b_{0}, \ldots, b_{d_{1} -1} , - b_{d_{1}} \right) $ is expected to be the shortest vector in $R_{k}$. The integers $b_{0}, \ldots, b_{d_{1}}$ give us a relations $\varphi_{2}$, such that for any tuple of coefficients $\left( x_{1}, \ldots, x_{n} \right)$, if $x_{1} = \theta $ is any root of $f_{1}$, then $x_{2} = \varphi_{2}\left( \theta \right)$. We repeat this with all coefficients. 

In fact, in all of our examples we found that for a given theta hyperplanes, all of the coefficients of its defining equation were defined over the number field generated by the minimal polynomial of one of the coefficients, so the method described above was sufficient to complete all of our computations. 

The correctness of the minimal polynomials and coefficient relations computed can be easily checked by verifying that the  points they define are solutions of the set of equations defining our zero-dimensional scheme $S$.

\subsection{Algebraic Hyperplanes and Complex Approximations}
\label{CA}
Suppose that $ \mathbf{a} = \left( a_{1}, \ldots, a_{n} \right)$ is a complex approximation, accurate to $k$ decimal places, of a point $\mathbf{x} = \left( x_{1}, \ldots , x_{n} \right)$ of $S$. In this subsection we give a method to compute the minimal polynomial of $x_{i}$. As in the $p$-adic case, we use the approximation to define a lattice and then search for short vectors in this lattice.
Computing the minimal polynomial of a complex approximation is implemented in \texttt{MAGMA} by the command \texttt{MinimalPolynomial}, whose input is a complex approximation and the supposed degree of its minimal polynomial, and its output is a polynomial with integer coefficients and of the given degree, and which has a root that is well approximated by our first input. This uses similar techniques to the method described in this subsection, but we found it to be slower, especially when the degree of the polynomial is large. For instance, our code was significantly faster when computing a degree $24$ minimal polynomial corresponding to the coefficient of one quadritangent to $X_{0}\left( 55 \right)$. We were able to run through all polynomials of degree $ 2 \le d \le 24$ in  $12.5$ seconds, whilst using the \texttt{MAGMA} command to complete the same task took $104.6$ seconds. This particular orbit of quadritangents was not used in our final calculation, but the file computing it is included in the online repository to demonstrate the efficiency of our implementation. 
In our examples, we chose generators of the $2-$torsion subgroup which correspond to theta hyperplanes whose coefficients had minimal polynomials of small degree, since this made our final computations and taking Galois invariants simpler. However, throughout our computations, many large degree polynomials were computed, from which optimal generators could be chosen. Thus, it was very important to have a program which computed the minimal polynomials as efficiently as possible. In all genus $5$ examples, we found that our implementation was at least three times faster, even when the degree of the polynomials was small. 

\subsubsection*{Minimal Polynomials}
Fix $i$, $1 \leq i \leq n$ and let $\theta = x_{i}$, $ a = a_{i}$. Recall that $a$ is a complex approximation of $\theta$, correct to $k$ decimal places. As $\theta$ is an algebraic number, there exists $d_{\theta} \in \mathbb{N}$ and $c_{0}, \ldots, c_{d_{\theta}} \in \mathbb{Z}$ such that 
\begin{center}
$c_{d_{\theta}}\theta^{d_{\theta}} + \ldots + c_{1}\theta + c_{0} = 0$.
\end{center}
If the imaginary part of $a$ is very small, it's likely that $a$ is approximating a real, algebraic number. In this case, we will take  $a  = \text{Real}\left( a \right) \in \mathbb{R}$ and let $C = 10^{k'}$, for some large natural number $k' < k$, such that  
\begin{center}
    $ \mid \left[ C\cdot a^{i} \right] - C \cdot \theta^{i} \mid \ \leq  \ 1 $ for all $0 \leq i \leq d_{\theta} $
\end{center}
 where $\left[ x \right]$ denotes the integer part of $x \in \R$. 
 
Let $\mathcal{L}$ be the lattice generated by the columns $v_{d_{\theta}}, \ldots, v_{0}$ of the matrix 
\begin{center}
$A = \begin{pmatrix}
1 & \ldots & 0 & 0  \\
0 & \ldots & 0 & 0 \\
\vdots & \ddots & \vdots & \vdots \\
0 & \ldots & 1& 0 \\
\left[ C a^{d_{\theta}} \right]   & \ldots & \left[ C a \right]  & \left[ C \right]  

\end{pmatrix}
 = \left( v_{d_{\theta}} , \ldots, v_{1}, v_{0} \right).
$
\end{center}
Let  $\gamma =   c_{d_{\theta}}\left[ C a^{d_{\theta}} \right]  + \ldots + c_{1} \left[ C a \right]  + c_{0}  \left[ C \right] $  and 
\begin{center}
$  \mathbf{c} = \begin{pmatrix}
c_{d_{\theta}} \\
\vdots \\
c_{1} \\
\gamma
\end{pmatrix} = c_{d_{\theta}}v_{n} + \ldots + c_{0}v_{0} .
$
\end{center}
Observe that $\mathbf{c} \in \mathcal{L}$ since $c_{d_{\theta}}, \ldots, c_{0} \in \Z$, and we can recover the vector of coefficients of the minimal polynomial $\mathbf{c}_{\infty} = \left( c_{d_{\theta}}, \ldots, c_{0} \right)^{T}$ from $\mathbf{c}$, by setting 
\begin{center}
$ c_{0} = \frac{1}{C} \left( \gamma - \left( c_{d_{\theta}}\left[ C a^{d_{\theta}} \right]  + \ldots + c_{1} \left[ C a \right]  \right) \right) $.
\end{center} 
Then 
\begin{center}
\begin{align*}
\vert \vert \mathbf{c} \vert \vert &  = \sqrt{ c_{d_{\theta}}^{2} + \ldots + c_{1}^{2} + \gamma^{2} } \\
&  \leq  \sqrt{ c_{d_{\theta}}^{2} + \ldots + c_{1}^{2} + \left( \gamma - Cc_{d_{\theta}}\theta^{d_{\theta}} - \ldots - Cc_{1}\theta - Cc_{0} \right)^{2} } \\
&  \leq  \sqrt{ c_{d_{\theta}}^{2} + \ldots + c_{1}^{2} + \left( c_{d_{\theta}}( \left[ Ca^{d_{\theta}} \right] - C\theta^{d_{\theta}})  + \ldots + c_{1} ( \left[ Ca \right] - C\theta ) + c_{0} ( \left[ C \right] - C) \right)^{2} } \\
&  \leq \sqrt{ c_{d_{\theta}}^{2} + \ldots + c_{1}^{2} + \left( c_{d_{\theta}} + \ldots + c_{1} + c_{0} \right)^{2} } \\
&  \leq \sqrt{ c_{d_{\theta}}^{2} + \ldots + c_{1}^{2} + \left( c_{d_{\theta}}^{2} + \ldots + c_{1}^{2} + c_{0}^{2} \right)^{2} } \\
&  \leq \sqrt{ 2 \left( c_{d_{\theta}}^{2} + \ldots + c_{1}^{2} + c_{0}^{2} \right)^{2} }  = \sqrt{2} \  \vert \vert \ \mathbf{c}_{\infty} \vert \vert^{2} 
\end{align*}
\end{center}
and this shows that although the length $\vert \vert \mathbf{c} \vert \vert$ depends on the precision of the approximation $k$,  $\vert \vert \mathbf{c} \vert \vert$ is bounded by the fixed constant $ \sqrt{2} \vert \vert \ \mathbf{c}_{\infty} \vert \vert^{2} $.

Hermite's theorem suggests that the length of the  shortest vector in $\mathcal{L}$ is approximately $\Delta \left( \mathcal{L} \right)^{\frac{1}{d_{\theta}+1 }}$, where $\Delta \left( \mathcal{L} \right)$ is the determinant of the lattice $\mathcal{L}_{k}$, as defined in \cite[Page 68]{smart}. As in the $p$-adic case, if $k$ is large enough we expect the vector $\mathbf{c}$ to be the shortest vector in the lattice, and of length significantly smaller than $\Delta \left( \mathcal{L} \right)^{\frac{1}{d_{\theta}+1 }}$. In our case, $\Delta \left( \mathcal{L} \right) = \text{det}\left( A \right) = C = 10^{k'}$; and so if our minimal polynomial has coefficients of order $10^{n}$, $k$, $k'$ are such that:
\begin{center}
 $\left( d_{\theta} + 1 \right)10^{2n} \leq 10^{k'/ \left( d_{\theta} + 1 \right)}$.
\end{center}
and if the shortest vector in $\mathcal{L}$ is shorter than $\Delta \left( \mathcal{L} \right)^{\frac{1}{d_{\theta}+1 }}$, then it is a suitable candidate for the vector we are looking for. As before, we search for the shortest vector in the lattice using the \texttt{Magma} command \texttt{ShortestVectors} (see \cite{MR1484478}).

\begin{rem}
When the imaginary part of $a$ is not small, the same method can be used but with $\mathcal{L}$ being generated by the columns of 
\begin{center}
$A_{k} = \begin{pmatrix}
1 & \ldots & 0 & 0   \\
0 & \ldots & 0 & 0   \\
\vdots & \ddots & \vdots & \vdots\\
0 & \ldots & 1 & 0 \\
\left[ C \text{Re}\left(a^{d_{\theta}} \right)  \right]   & \ldots  &   \left[ C  \text{Re} \left( a \right)  \right] &  \left[ C \right]  \\
\left[ C \text{Im}\left(  a^{d_{\theta}} \right)  \right]   & \ldots   &   \left[ C  \text{Im} \left( a \right)  \right] &  0
\end{pmatrix}
$
\end{center}
where $\text{Re} \left( x \right)$ and $\text{Im} \left( x \right)$ denote the real and imaginary parts of $x \in \mathbb{C}$; and $C = 10^{k'}$ where $k^{'} < k $, is such that 
\begin{center}
      $ \mid \left[ C\cdot \text{Re}\left( a^{i}\right) \right] - C \cdot \text{Re} \left( \theta^{i} \right)  \mid \ \leq  \ 1 $  and  $ \mid \left[ C\cdot \text{Im}\left( a^{i}\right) \right] - C \cdot \text{Im} \left( \theta^{i} \right)  \mid \ \leq  \ 1 $       
      \end{center}
for all $0 \leq i \leq d_{\theta} $. Arguing as in the previous case, for suitable $k$ and $k^{'}$, the shortest vector in $\mathcal{L}$ is a suitable candidate for vector corresponding to the minimal polynomial.
\end{rem}

As in the $p$-adic case, we can replace $d_{\theta}$ in the above expression by any positive integer $d$ and search for short vectors in the lattice constructed  with respect to this $d$. The strategy to find $\mathbf{c}$ is as before, we choose a candidate $d$ for $d_{\theta}$, starting with $d=1$, and running through the natural numbers. If our candidate $d$ is equal to $d_{\theta}$, then the shortest vector and the polynomial whose coefficients are derived from the shortest vector should satisfy the same conditions as before. 

\subsubsection*{Coefficient Relations}
As in the previous subsection, once we have candidates $ \left( f_{1}, \ldots, f_{n} \right)$ for the minimal polynomials of $\left( x_{1}, \ldots, x_{n} \right)$, we want to find  which tuples of roots correspond to coefficients of equations defining theta hyperplanes. 

Factorising in this case is simpler, as we now factor $f_{i}$  over $\mathbb{C}$,
\begin{center}
$f_{i} = s_{1} \ldots s_{d_{i}}$.
\end{center}
Then for some $n$, $s_{n}\left( a_{i} \right)$ is very small. Thus $s_{n}$ corresponds to the root of $f_{i}$ approximated by $a_{i}$. 

 When factorisation is not efficient, we search for relations amongst coefficients, say between $a_{1}$ and $a_{2}$ using lattice reductions. If $a_{1}, a_{2} \in \mathbb{R}$, that is the imaginary part of both approximations is very small, we search for $b_{0}, \ldots, b_{d_{1}} $ such that 
\begin{center}
    $-b_{d_{1}}x_{2} = b_{d_{1} -1}x_{1}^{d_{1} -1} + \ldots + b_{1}x_{1} + b_{0}$
\end{center}
where $d_{1}$ is the minimal polynomial of $f_{1}$. Similar to the method used to search for minimal polynomial, such integers can be found by searching for short vectors in the lattice generate by the columns of 
\begin{center}
$A = \begin{pmatrix}
1 & \ldots & 0  & 0   & 0 \\
0 & \ldots & 0 & 0 & 0  \\
\vdots & \ddots & \vdots & \vdots& \vdots  \\
0 & \ldots &  0 & 1 &  0 \\
\left[ C a_{1}^{d_{1}-1} \right]   & \ldots & \left[ C a_{1}  \right]   & \left[ C a_{2}  \right] & \left[ C \right]  

\end{pmatrix}
$
\end{center}
where $C = 10^{k'}$ and $k' < k$ is chosen as before. If the imaginary parts of $a_{1}, a_{2} $ are not both small, we instead search for short vectors in the lattice generated by the columns of 

\begin{center}
$A = \begin{pmatrix}
1 & \ldots & 0 & 0 & 0   \\
0 & \ldots & 0 & 0 & 0  \\
\vdots & \ddots & \vdots & \vdots & \vdots \\
0 & \ldots & 0 & 1 & 0 \\
\left[ C \text{Re}\left(  a_{1}^{d_{1}-1} \right)  \right]   & \ldots & \left[ C  \text{Re} \left( a_{1} \right)  \right]  &   \left[ C  \text{Re} \left( a_{2}  \right)  \right] &  \left[ C \right]  \\
\left[ C \text{Im}\left(  a_{1}^{d_{1}-1} \right)  \right]   & \ldots & \left[ C  \text{Im} \left( a_{1} \right)  \right]  &   \left[ C  \text{Im} \left( a_{2} \right)  \right] & 0
\end{pmatrix}
$
\end{center}

In all of our examples, for a given orbit of theta hyperplanes, we were able to express all coefficients of the defining equations in terms of one fixed coefficient. 

The correctness of the minimal polynomials and coefficient relations computed can be verified by computing the corresponding points and checking that they are solutions to the set of equations defining our zero-dimensional scheme.

\section{Examples}
\allowdisplaybreaks
\subsection{The genus $3$ curve: $X_{0} \left( 75 \right) / w_{25}$}
\label{g3eg}
Consider the non-hyperelliptic genus $3$ curve $X$ over $\Q$ defined by:
\begin{equation*}
 f( x,y,z) = 3x^3z - 3x^2y^2 + 5x^2z^2 -3xy^3 -19xy^2z -xyz^2 + 3xz^3+ 2y^4 + 7y^3z -7y^2z^2 -3yz^3.
\end{equation*}
This is the quotient of the  modular curve $X_{0} ( 75 )$ by the action of the Atkin-Lehner operator $w_{25}$ computed in \cite[Pages 19-21]{freitas2015elliptic}.
Let $J$ be the Jacobian variety of $X$. In \cite{freitas2015elliptic} the authors showed $J\left(\mathbb{Q} \right)_{\text{tors}} \cong  \ \mathbb{Z}/2\mathbb{Z} \times \mathbb{Z} / 40\mathbb{Z}$. Using the methods described in Sections 3,4 and 5 we show $J\left( \mathbb{Q} \right) \left[ 2 \right] \cong  \left( \mathbb{Z} / 2 \mathbb{Z} \right)^{2}$.

First, we write down equations for the scheme of bitangents. We work on the affine chart $ \lbrace z =1 \rbrace$ and with bitangents of the form $x = a_{1}y + a_{2}$ for some $a_{1},a_{2} \in \overline{\mathbb{Q}}$. Using the method of Section 3.1, we obtain the equations:

\begin{align*}
 e_{1}  & = 3a_{1}^{2}a_{4}^{2} + 3a_{1}a_{4}^{2} + 3a_{2}^3 + 5a_{2}^{2} + 3a_{2} - 2a_{4}^{2}, \\
e_{2} &= 6a_{1}^{2}a_{3}a_{4} + 9a_{1}a_{2}^{2} + 10a_{1}a_{2} + 6a_{1}a_{3}a_{4} +  3a_{1} - a_{2} - 4a_{3}a_{4} - 3,\\
e_{3}  & = 9a_{1}^{2}a_{2} + 3a_{1}^{2}a_{3}^{2} + 6a_{1}^{2}a_{4} + 5a_{1}^{2} + 3a_{1}a_{3}^{2} +  6a_{1}a_{4} - a_{1} - 3a_{2}^{2} - 19a_{2} - 2a_{3}^{2} - 4a_{4} - 7,\\
e_{4}  &= 3a_{1}^{3} + 6a_{1}^{2}a_{3} - 6a_{1}a_{2} + 6a_{1}a_{3} - 19a_{1} - 3a_{2} -  4a_{3} + 7.
\end{align*}

Let $S$ be the scheme defined by $e_{1},e_{2},e_{3},e_{4}$. As $17$ is a prime of good reduction for $X$ and we can verify using \texttt{Magma} that $S /  \mathbb{F}_{289}$ has degree 28; all bitangents to $X$ correspond to points on $S$.

Let $L = \mathbb{Q} \left( \sqrt{7} \right)$ and $\mathcal{O}_{L} = \mathbb{Z} \left[ \sqrt{7} \right]$ its ring of integers. The prime ideal $\mathfrak{p} =  \langle 17 \rangle $ of $\mathcal{O}_{L}$ has norm $17^{2} = 289$, so we can use the isomorphism $\mathbb{F}_{289} \cong \mathcal{O}_{L}/\mathfrak{p} = \mathbb{F}_{\mathfrak{p}}$ to calculate the points of $S \left( \mathbb{F}_{289}\right) =  S\left( \mathbb{F}_{\mathfrak{p}} \right)$.

Each point in $S\left( \mathbb{F}_{\mathfrak{p}} \right)$ has a unique lift modulo $\mathfrak{p}^{k}$ for each $k \geq 1$ and using the methods described in Section 5.1, we find the minimal polynomials of the coefficients of each bitangent, and define the bitangents. 

There are four bitangents defined over the rationals, given by the following equations 
\begin{align*}
 & x + 3y + z, \\
 & x - y + z, \\
& x - 2y +  z, \\
& x + 2y +z.  
\end{align*}
There are six bitangents defined over the quadratic fields $\Q \left(\sqrt{5}\right)$, with equations:
\begin{align*}
    & -2x + \left( -1 \pm \sqrt{5} \right) y + \left( 1 \pm 3 \sqrt{5} \right) z, \\
    & -2x + \left( 7 \pm \sqrt{5} \right) y + \left( -3 \mp \sqrt{5} \right) z, \\
    & x + \left( -4 \pm 2 \sqrt{5} \right) y + \left( -9 \pm 4 \sqrt{5} \right) z.
\end{align*}
The remaining 18 bitangents are listed in Galois orbits in Table \ref{table: 1}. All bitangents are described by equations of the form  $x - a_{1} y - a_{2}z$ and we give a minimal polynomial $f$ of $a_{1}$ and an expression for $a_{2}$ in terms of $a_{1}$. 
\allowdisplaybreaks
\begin{center}
\begin{table}
\caption{Equations defining some bitangents of $X_{0}\left( 75 \right)/ w_{25}$.} 
\label{table: 1}
\begin{align*}
\hline
&  f = u^6 + 3u^5 - 5u^3 + 60u^2 + 63u - 41 \\
\hline 
& a_{1} = u  \\
& a_{2} = \frac{1}{27}\left(  u^5 + u^4 - 2u^3 - u^2 + 62u  - 34 \right)  \\
\hline 
\hline
& f = 11u^6 + 49u^5 + 20u^4 - 55u^3 + 40u^2 - 11u + 1 \\
\hline 
& a_{1} = u  \\
& a_{2} = \frac{1}{17}\left( 407u^5 + 1901u^4 + 1154u^3 - 1773u^2 + 1106u  - 188 \right) \\
\hline 
\hline
& f = u^{6} + 14u^{5} - 5u^{4} - 20u^{3} + 95u^{2} + 134u - 139 \\
\hline 
 & a_{1} = u  \\
& a_{2} = \frac{1}{9164}\left( 41u^5 + 467u^4 - 1312u^3 + 2604u^2 + 2687u - 3083 \right) \\
\hline 
\end{align*}
\end{table}
\end{center}

The minimal polynomials show that the Galois orbits of the bitangents are as follows:
\begin{itemize}
\item 3 orbits with 6 bitangents each;
\item 3 orbits with 2 bitangents each, all defined over $\Q \left( \sqrt{5} \right)$;
\item 4 bitangents defined over $\mathbb{Q}$ and thus stable under Galois action.
\end{itemize}
Let $K$ be the splitting field of  $u^{6} + 14u^{5} - 5u^{4} - 20u^{3} + 95u^{2} + 134u - 139$. This is a degree 12 Galois extension of $\mathbb{Q}$ and we verify using the \texttt{Magma} command \texttt{Roots}, see \cite{MR1484478}, that all of the coefficients stated in Table \ref{table: 1} split over this number field and thus all bitangents are defined over $K$. Furthermore, the same command can be used to verify that this is the smallest degree number field over which all bitangents are defined. Let $\mathcal{O}_{K}$  be the ring of integers of $K$. We view $X$ as a projective curve over $K$ and scale the equations of the bitangents to ensure that they are defined over the maximal order of $K$. Denote the set equations, defined over $\mathcal{O}_{K}$, which cut out the bitangents to $X$ by BT. 

Fix the rational bitangent 
\begin{center}
$b : x + 3y + z = 0 $
\end{center}
and define $HBT = \lbrace \frac{1}{2} \text{div} \left( \frac{l}{b} \right) : l \in BT \rbrace $. Let $H$ be the subgroup of $J \left( K \right)  \left[ 2 \right]$ generated by equivalence classes of elements of $HBT$.

The ideal $\mathfrak{p} = \langle 17, 4 + 2\theta + 2 \theta^2 \rangle$ is a prime ideal in $ \mathcal{O}_{K}$ of norm $289$, where $\theta$ is a generator of $K$.  As $17$ is a prime of good reduction for $X$, reduction modulo $\mathfrak{p}$ induces an injection 
\begin{center}
$r_{\mathfrak{p}} : J \left( K \right) _{\text{tor}} \longrightarrow J \left( \mathbb{F_{\mathfrak{p}}} \right)$
\end{center}
see \cite{katz}. It can be shown using \texttt{Magma} (see the examples of \cite[Section 12.9]{BPS} for similar computation or refer to the accompanying code for details) that $ J \left( \mathbb{F_{\mathfrak{p}}} \right) \cong \left( \mathbb{Z} / 2 \mathbb{Z} \right) \times \left( \mathbb{Z} /8 \mathbb{Z} \right)^{2} \times \left( \mathbb{Z} / 40 \mathbb{Z} \right)^{2} \times \left( \mathbb{Z} / 160 \mathbb{Z} \right)$ and $r_{\mathfrak{p}} \left( H \right) \cong \left( \mathbb{Z} / 2 \mathbb{Z} \right)^{6}$, and thus we conclude that $H = J \left( K \right)\left[ 2 \right]$. As the genus of $X$ is 3, $H$ is necessarily the entire $2$-torsion subgroup,
\begin{center}
$H = J \left( \overline{\mathbb{Q}} \right) \left[ 2 \right]$.
\end{center}
To determine the rational $2$-torsion subgroup, we take Galois invariants of $H$. Let $G = \Gal \left( K/ \mathbb{Q} \right)$.  Elementary computations show $G  \cong D_{12} $, where $D_{12}$ is the dihedral group of order 12, and we find that the subgroup of $H$ fixed by $G$ is 
\begin{center}
$H^{G} =  \left( \mathbb{Z} / 2 \mathbb{Z} \right)  \left[ P_{1} + P_{2} -2P_{0} \right]  \ + \   \left( \mathbb{Z} / 2 \mathbb{Z} \right)  \left[ P_{3} + P_{4} - 2P_{0} \right] $
\end{center}
where the points $P_{0},\ldots, P_{4}$ are the points of the intersection of the first $3$ rational bitangents stated earlier in this subsection. 
\begin{align*}
 & P_{0} = \left( 1 :-1 : 2 \right), \\
 & P_{1} = \left( 1 - \sqrt{2} :1 : \sqrt{2}\right), \quad  P_{2} =\left( 1 + \sqrt{2} :1 : - \sqrt{2}\right), \\
 & P_{3} = \left( -6 + 2 \sqrt{-15} : 1 + \sqrt{-15} : 8 \right), \quad  P_{4} = \left( -6 - 2 \sqrt{-15} : 1 - \sqrt{-15} : 8 \right).
\end{align*}

Thus $J \left( \mathbb{Q} \right) \left[ 2 \right] \cong \left( \mathbb{Z} / 2 \mathbb{Z} \right)^{2}$, verifying the result of \cite{freitas2015elliptic}.

\subsection{The genus $4$ curve: $X_{0}\left( 54 \right)$}
In this subsection we compute the $2$-torsion subgroup of the Jacobian of the non-hyperelliptic genus 4 modular curve $X = X_{0}\left( 54 \right)$. This has a model over the rationals given by:
\begin{align*}
f_{1} & = x^{2}z - xz^{2} - y^{3} + y^{2}w -3yw^{2} + z^{3} + 3w^{3}, \\
f_{2} &= xw - yz + zw
\end{align*}
as given in \cite{OS}). Let $J$ be the Jacobian of this curve. As stated in \cite{OS},  $J \left(  \mathbb{Q }\right) \cong \left( \mathbb{Z} / 3 \mathbb{Z} \right) \times \left( \mathbb{Z} / 3 \mathbb{Z} \right) \times \left( \mathbb{Z} / 9 \mathbb{Z} \right) $, and so there is no rational 2-torsion. We verify this using our method. 

For the scheme of tritangents, we work on the affine chart $\lbrace w = 1 \rbrace$ and with tritangent planes of the form 
\begin{center}
$x = a_{1}y + a_{2} z + a_{3}$
\end{center}
for some $a_{1}, a_{2}, a_{3} \in \mathbb{\overline{Q}}$. Using the method described in Section 3.2 we define a zero-dimensional scheme of tritangent planes to $X$, which we call $S$. It can be verified using \texttt{Magma} that $S$ is non-singular over $\mathbb{F}_{289}$ and it has degree $120$, thus all tritangents to $X$ correspond to a point on $S$.

Let $L = \mathbb{Q} \left( \sqrt{-3} \right)$ and $\mathcal{O}_{L}$ its ring of integers. The ideal $ \mathfrak{p} = \left\langle 17 \right\rangle $ is prime and its norm is 289,  and we use the isomorphism $ \mathbb{F}_{289} \cong \mathbb{F}_{\mathfrak{p}} = \mathcal{O}_{L} / \mathfrak{p}$ to calculate the points of $ S\left ( \mathbb{F}_{289} \right) = S \left( \mathbb{F}_{\mathfrak{p} } \right)$. 

Using the techniques of Sections 4 and 5, we compute lifts of the points of $S\left( \mathbb{F}_{\mathfrak{p}} \right)$ modulo powers of $\mathfrak{p}$, and compute the minimal polynomials.  In this example, computing the minimal polynomials of all $a_{1}$s was sufficient to determine the field of definition of the $2$-torsion subgroup and to complete our calculation. 
The minimal polynomials show that the Galois orbits of tritangents are as follows:
\begin{itemize}
\item there are 2 orbits with 3 tritangents each;
\item there is one orbit with 6 tritangents;
\item there are 2 orbits with 9 tritangents each; 
\item there are 3 orbits with 18 tritangents each; 
\item there is 1 orbits with 36 tritangents.
\end{itemize}
 
 The number field over which all of our minimal polynomials split is a candidate for the field $K = \Q \left( J \left[ 2 \right] \right)$. This has degree $36$ and we were able to find all $120$ points of $S$ over $K$ using the \texttt{Magma} command \texttt{Points}, confirming that our candidate is $K$. From this we deduce the $2$-torsion subgroup $J\left[ 2 \right]$ as before.

Let $G = \Gal \left( K / \mathbb{Q} \right)$. Elementary computation show that $G \cong C_{3} \times A_{4}$, where $C_{3}$ denotes the cyclic group of order $3$ and $A_{4}$ denotes the alternating group acting on $4$ elements. Let $s \in G $ be any element of order $6$. We find that $s$ has no fixed points so $J \left[ 2 \right] ^{s} = 0$. Therefore, $J\left( \mathbb{Q} \right) \left[ 2 \right] = 0 $, verifying the result of \cite{OS}. 

\subsection{The genus $5$ curve: $X_{0}\left( 42 \right)$}
\label{g5ex}
In this subsection we compute the $2$-torsion subgroups of the Jacobian of the genus $5$, non-hyperelliptic modular curve $X = X_{0}\left( 42 \right)$. The canonical model that we use is (as in \cite{OS}) the intersection of the three quadrics:
\begin{align*}
f_{1} & = x_{1}x_{3} - x_{2}^{2} + x_{3}x_{4}; \\
f_{2} & = x_{1}x_{5} - x_{2}x_{5} - x_{3}^{2} + x_{4}x_{5} - x_{5}^{2};  \\
f_{3} & = x_{1}x_{4} - x_{2}x_{3} + x_{2}x_{4} - x_{3}^{2} + x_{3}x_{4} + x_{3}x_{5} - x_{4}^{2} - 2x_{4}x_{5}.
\end{align*}
For the scheme of quadritangents,  with this model, it is convenient to work on the affine chart $ \lbrace  x_{5} = 1 \rbrace$  and with planes of the form 
\begin{center}
$x_{1} = a_{1}x_{2} + a_{2}x_{3} + a_{3}x_{4} + a_{4}$ 
\end{center}
for some $a_{1},a_{2}, a_{3},a_{4} \in \overline{\mathbb{Q}}$.
The intersection of such a plane with the affine curve is given by 
\begin{center}
$F_{i}\left( x_{2}, x_{3}, x_{4} \right) = f_{i} \left( a_{1}x_{2} + a_{2}x_{3} + a_{3}x_{4} + a_{4}, x_{2} , x_{3}, x_{4}, 1 \right) $  for $i = 1,2,3 .$
\end{center}
Explicitly, these expressions are 
\begin{center}
$F_{1} \left(x_{2}, x_{3}, x_{4} \right) = -x_{2}^{2} + a_{1} x_{2}x_{3} + a_{2}x_{3}^{2} + \left( a_{3} + 1 \right) x_{3}x_{4} + a_{4}x_{3}; $\\
$F_{2} \left( x_{2} , x_{3} , x_{4} \right) = \left( a_{1} - 1\right) x_{2} -x_{3}^{2} + a_{2}x_{3} + \left( a_{3} + 1 \right) x_{4} + a_{4} -1; $\\
$F_{3} \left( x_{2}, x_{3}, x_{4} \right) = -x_{2}x_{3} + \left( a_{1} + 1 \right) x_{2}x_{4} -x_{3}^{2} + \left( a_{2} + 1 \right) x_{3}x_{4} + x_{3} + \left( a_{3} - 1 \right) x_{4}^{2} + \left( a_{4} - 2 \right) x_{4}. $
\end{center}
For $a_{1} \neq 1 $, $F_{2} = 0 $ gives an expression for $x_{2}$ in terms of $x_{3}$ and $x_{4} $
\begin{align}
x_{2}   = \frac{-1}{a_{1} - 1 } \left( -x_{3}^{2} + a_{2} x_{3} + \left( a_{3} + 1 \right) x_{4} + a_{4} - 1 \right)   
   = \phi \left( x_{3} , x_{4} \right).
  \end{align}
  
 Substituting for $x_{2} $ in $F_{1}$ and $F_{3}$, we get 
 \begin{center}
 \begin{align*}
 \widetilde{G_{1}} & = F_{1} \left( \phi \left( x_{3}, x_{4} \right), x_{3}, x_{4} \right) 
  =  \frac{G_{1} \left( x_{3}, x_{4} \right) }{\left( a_{1} - 1 \right)^{2}};  \\
   \widetilde{G_{2}} & = F_{3} \left( \phi \left( x_{3}, x_{4} \right), x_{3}, x_{4} \right) 
  =  \frac{G_{2} \left( x_{3}, x_{4} \right) }{\left( a_{1} - 1 \right)}.
 \end{align*}
 \end{center}

Clearing denominators, the intersection is given by 
\begin{center}
$ G_{1} \left( x_{3}, x_{4} \right) = G_{2} \left( x_{3} , x_{4} \right) =0  $ 
\end{center}
where 
 \begin{center} 
 
$ G_{1} = -x_{3}^{4} + \left(a_{1}^{2} - a_{1}+ 2a_{2}\right)x_{3}^{3} + \left(2a_{3} + 2\right)x_{3}^{2}x_{4} + \left(-a_{1}a_{2} - a_{2}^{2} + a_{2} + 2a_{4} - 2\right)x_{3}^{2} + \left(-a_{1}a_{3} - a_{1} -  2a_{2}a_{3} - 2a_{2} + a_{3} + 1\right)x_{3}x_{4} + \left(a_{1}^{2} - a_{1}a_{4} - a_{1} -  2a_{2}a_{4} + 2a_{2} + a_{4}\right)x_{3} + \left(-a_{3}^{2} - 2a_{3} - 1\right)x_{4}^{2} +  \left(-2a_{3}a_{4} + 2a_{3} - 2a_{4} + 2\right)x_{4} - a_{4}^{2} + 2a_{4} - 1;$ \\
$G_{2} = -x_{3}^{3} + \left(a_{1} + 1\right)x_{3}^{2}x_{4} + \left(-a_{1} + a_{2} + 1\right)x_{3}^{2} + \left(a_{1} - 2a_{2} + a_{3}\right)x_{3}x_{4} + \left(a_{1} + a_{4} - 2\right)x_{3} + \left(-2a_{1} - 2a_{3}\right)x_{4}^{2} + 
   \left(-a_{1} - 2a_{4} + 3\right)x_{4}.$
\end{center}
These can be re-written as 
\begin{center}
$G_{1} = g_{1} \left( x_{3} \right) + h_{1} \left( x_{3} \right) x_{4} + \alpha_{1} x_{4}^{2}; $ \\
$G_{2} =  g_{2} \left( x_{3} \right) + h_{2} \left( x_{3} \right) x_{4} + \alpha_{2} x_{4}^{2} $
\end{center}
with $ g_{1}, g_{2}, h_{1}, h_{2} \in \mathbb{Q} \left[ a_{1}, a_{2}, a_{3} ,a_{4} \right] \left[ x_{3} \right] $ and $ \alpha_{1} = -\left( a_{3} + 1 \right)^{2} $,  $ \alpha_{2} = -2\left( a_{1} + a_{2} \right) $. 
If $\alpha_{1} \neq 0 $ and $ \alpha_{2} \neq 0 $ then
\begin{center}
$ \alpha_{2} G_{2} - \alpha_{1} G_{2} =  \left( \alpha_{2} h_{1} \left( x_{3} \right)  - \alpha_{1} h_{2} \left( x_{3} \right) \right) x_{4} + \alpha_{2} g_{1} \left( x_{3} \right) - \alpha_{1} g_{2} \left( x_{3} \right) = 0 $ 
\end{center}
and if $ T_{1}  =  \alpha_{2} h_{1} \left( x_{3} \right)  - \alpha_{1} h_{2} \left( x_{3} \right) \neq 0  $ , as a polynomial in $x_{3} $,  the above gives an expression for $x_{4}  = \frac{-T_{2} }{T_{1} }$ with 
$ T_{2} =    \alpha_{2} g_{1} \left( x_{3} \right)  -  \alpha_{1} g_{2} \left( x_{3} \right) $.
Substituing for $x_{4} = -T_{2}/T_{1}$ in $G_{1}$, we get 
\begin{center}
$G\left( x_{3} \right) = G_{1} \left( x_{3} , \frac{-T_{2}}{T_{1} } \right) =  \frac{h\left( x_{3} \right) }{g \left( x_{3} \right) } $
\end{center}
where $h \in \mathbb{Q}\left[ a_{1}, a_{2} , a_{3}, a_{4} \right]\left[ x_{3} \right]$ has degree 8 and $g \left( x_{3} \right) =  T_{1}^{2}/ \left(a_{3} + 1  \right)^{2} $. 

Clearing denominators in the expression above,  we remark that if the given plane $x_{1} = a_{1}x_{2} + a_{2}x_{3} + a_{3}x_{4} + a_{4} $ is a quadritangent,  then $h$ is necessarily a square.  Equivalently,  there exist $a_{5} , a_{6}, a_{7}, a_{8} \in \overline{\mathbb{Q}}$ such that 
\begin{center}
$h\left( x \right) = l \left( x^{4} + a_{5}x^{3} + a_{6}x^{2} + a_{7}x + a_{8} \right)^{2} $
\end{center}
where $l$ is the leading coefficient of $h$.  Equating coefficients in the above expression, gives 8 equations $e_{1} , \ldots ,e_{8} $ in $a_{1}, \ldots , a_{8}$.

We also add a 9th equation (and a 9th variable)
\begin{center}
$e_{9} : a_{9} \Delta \left( x^{4} + a_{5} x^{3} + a_{6}x^{2} + a_{7}x + a_{8} \right) + 1  = 0 $ 
\end{center}
to ensure that $ x^{4} + a_{5} x^{3} + a_{6}x^{2} + a_{7}x + a_{8}$ has non-zero discriminant and avoid singularities on our scheme. 
To derive $e_{1}, \ldots, e_{8}$, we assumed that $ a_{1} \neq 1$, $\alpha_{1} \neq 0 $ and $\alpha_{2} \neq 0 $,  and equations are also required for these conditions:
\begin{align*}
e_{10} & : a_{10} \left( a_{1} - 1 \right) + 1 = 0,   \\
e_{11} & : a_{11} \left( a_{3} + 1 \right) + 1 = 0, \\
e_{12} & : a_{12} \left( a_{1} + a_{3} \right) + 1 = 0. 
\end{align*}
It can be checked that $e_{2} \ldots e_{12}$ are irreducible, and $e_{1} = s_{1}s_{2} $, 
\begin{center}
$s_{1} = a_{3}a_{4} -a_{3}a_{8} - a_{3} -2a_{4}^{2} + 5a_{4} + a_{8} -3, $ \\
$s_{2} = a_{3} a_{4} + a_{3} a_{8} - a_{3} - 2a_{4}^{2} + 5a_{4} - a_{8} - 3. $
\end{center}
\texttt{Julia} suggests that the system:
\begin{center}
$s_{1} = s_{2} = e_{2} = \ldots = e_{12} = 0 $  
\end{center}
has no solutions.  We consider the cases $s_{1} =0 $ and $s_{2} =0 $ separately. 

The condition $T_{1} \neq 0$ will also require equations. As a polynomial, $T_{1} $ is non-zero, if its coefficients $t_{1}$, $t_{2}$ and $t_{3} $ are not all zero. 
\begin{align*}
t_{1} & = a_{1} a_{3} - 3a_{1} -3a_{3} + 1, \\
t_{2} &= k_{1} k_{2} = \left( a_{1} + 2a_{2} + a_{3} \right) \left( 2a_{1} + a_{3} - 1 \right),\\
t_{3} &= a_{1}a_{3} - 4a_{1} a_{4} + 5a_{1} -2a_{3}a_{4} + a_{3} + 2a_{4} - 3.
\end{align*}
\subsubsection*{Case 1}
The first 12 equations are  $s_{1}, e_{2}, \ldots e_{12} $ and we consider all possible combinations of zero and non-zero $t_{1}, k_{1}, k_{2} $ and $t_{3}$.

\textbf{Case 1.1 :  $t_{1} \neq 0, k_{1} \neq 0, k_{2} \neq 0, t_{3} \neq 0 $ }
We add 4 equations and 4 variables 
\begin{align*}
e_{13} & : a_{13}t_{1} + 1 = 0; \\
e_{14} & : a_{14}k_{1} + 1 = 0; \\
e_{15} & : a_{15}k_{2} + 1 = 0;\\
e_{16} & : a_{16}t_{3} + 1 =0.
\end{align*}
\texttt{Julia} finds 96 approximate solutions to the system formed by the 16 equations in $ a_{1} , \ldots , a_{16} $.

\textbf{Case 1.2 :  $t_{1} \neq 0, k_{1} \neq 0, k_{2} =  0, t_{3} \neq 0 $ }
We add 4 equations and 3 variables 
\begin{align*}
e_{13} & : a_{13}t_{1} + 1 = 0; \\
e_{14} & : a_{14}k_{1} + 1 = 0; \\
e_{15} & : a_{15}t_{3} + 1 = 0; \\
e_{16} & : k_{2}  = 0.
\end{align*}
\texttt{Julia} finds 24 approximate solutions to the  system formed by the 16 equations in $ a_{1} , \ldots , a_{15} $.

\textbf{Case 1.3 :  $t_{1} \neq 0, k_{1}  =  0, k_{2} \neq  0, t_{3} \neq 0 $ }
We add 4 equations and 3 variables 
\begin{align*}
e_{13} & : a_{13}t_{1} + 1 = 0; \\
e_{14} & : a_{14}k_{2} + 1 = 0; \\
e_{15} & : a_{15}t_{3} + 1 = 0; \\
e_{16} & : k_{1}  = 0.
\end{align*}
\texttt{Julia} finds 14 approximate solutions to the system formed by the 16 equations in $ a_{1} , \ldots , a_{15} $.

In all other cases,  for all other possible combinations of zero and non-zero $t_{1}, k_{1}, k_{2}, t_{3} $, the corresponding systems have no approximate solutions. 

\subsubsection*{Case 2 }
The first 12 equations are $s_{2}, e_{2} , \ldots,e_{12} $ . 

\textbf{Case 2.1 :  $t_{1} \neq 0, k_{1} \neq 0, k_{2} \neq  0, t_{3} \neq 0 $ }
We add 4 equations and 4 variables 
\begin{align*}
e_{13} & : a_{13}t_{1} + 1 = 0; \\
e_{14} & : a_{14}k_{1} + 1 = 0; \\
e_{15} & : a_{15}k_{2} + 1 = 0; \\
e_{16} & : a_{16}t_{3} + 1 = 0. 
\end{align*}
\texttt{Julia} finds 256 approximate solutions to the  system formed by the 16 equations in $ a_{1} , \ldots , a_{16} $.

In all other cases,  the resulting systems have no approximate solutions.  We can also derive equations and schemes in the extreme cases, $\alpha_{1} =0$, $\alpha_{2} = 0$,  $a_{1} = 1$ etc. These cases combined had few solutions, and in fact these planes are not needed in our calculation of the 2-torsion subgroup. 

In each case,  we approximate solutions as complex points and using the techniques described in Section 5.2,  we compute the corresponding quadritangents to $X$. The 16 planes described in Table \ref{table: 3} and the cusps of $X_{0}\left( 42\right)$ were sufficient to fully describe the 2-torsion subgroup of $X_{0} \left( 42 \right)$. The planes computed occur in two Galois orbits. For both orbits, we give the minimal polynomial of $a_{1}$ and expression for $a_{2}, a_{3}, a_{4}$ in terms of $a_{1}$.

\begin{center}
\begin{table}
\caption{Some Quadritangent Planes to $X_{0}\left( 42 \right)$} 
\label{table: 3}
\begin{align*}
\hline & u^8 + 14u^7 + 151u^6 - 396u^5 + 283u^4 - 1730u^3 + 3201u^2 - 1440u +2284 \\
 \hline 
 a_{1} & = u \\
a_{2} & =  1/301100972656( 39380331u^7 + 442881623u^6 + 4027975134u^5 - 37845583334u^4  \\ & - 
    12023416509u^3    - 80118862717u^2 - 49248922084u + 241405142628 ) \\
 a_{3} & =  -1 /  301100972656 ( 52272255u^7 + 716826167u^6 + 7618567014u^5 - 24269489250u^4 \\ & + 8461415807u^3    - 144514109641u^2 + 341498511372u - 30420302668 ) \\
 a_{4} & =  -1/301100972656 (16239354u^7 + 285678727u^6 + 3306671724u^5 + 3184311718u^4  \\ & - 7648219606u^3 - 17837849633u^2 + 78572314624u  - 372802926204 ) \\
 \hline \\
 \hfill \\
 \hline 
 & 23u^8 + 78u^7 + 135u^6 + 146u^5 + 236u^4 + 322u^3 + 239u^2 + 94u + 23 \\
 \hline 
 a_{1} & = u \\
 a_{2} & = 1/6324 ( 10258u^7 + 27083u^6 + 40152u^5 + 34112u^4 + 75578u^3 + 80920u^2 \\ & + 41542u + 5979) \\
 a_{3} & =  -1/4743 (4324u^7 + 9535u^6 + 12034u^5 + 6080u^4 + 22416u^3 + 14876u^2 \\ & + 3588u - 7468) \\
 a_{4} & = -1/9486 ( 3404u^7 + 10601u^6 + 16874u^5 + 17098u^4 + 28680u^3 + 37276u^2 \\ &  + 24084u  - 7931) \\
 \hline
\end{align*}
\end{table}
\end{center}

Let $K$ be the number field defined by $f$
\begin{center}
$
 f\left( u \right) = 2713u^{16} +   9264u^{15} + 24252u^{14} - 1352u^{13} - 270446u^{12} - 739224u^{11} - 599968u^{10} + 1502208u^{9} + 
 6136803u^{8} + 10670696u^{7} + 11231488u^{6} + 7603968u^{5} + 3052898u^{4} + 591416u^{3} + 141500u^{2} + 210760u + 154057 
$.
\end{center}
This is a degree 16 Galois extension, and the above quadritangents are all defined over $K$. Let $\mathcal{O}_{K}$ be the ring of integers of $K$. The equations of the planes described above can be homogenised and scaled appropriately to ensure that they are all defined over $\mathcal{O}_{K}$. Denote by $QT$ the set of  equations of 16 quadritangents to $X$. The quotient of any $2$ elements of $QT$ is a function on the curve $X$. Fix an element $l_{0} \in QT$ and let 
\begin{center}
$HQT = \lbrace \frac{1}{2} \text{div} \left( \frac{l}{l_{0}} \right) : l \in QT \rbrace $.
\end{center}

The cusps of  $X_{0}\left( 42 \right) $ are:
\begin{align*}
c_{1}  & =  \left( 2 : -1 : 1 : -1 : 1 \right),  \\
c_{2} & =  \left( 2 : 2 : 1 : 2: 1 \right), \\
c_{3} & =  \left( 7 : 4 : 2 : 1 : 2 \right), \\ 
c_{4} & =  \left( 1 : -2 : 2 : 1 : 2 \right), \\  
 c_{5} & =  \left( 3 : 0 : 0 : -1 : 2 \right), \\
 c_{6} & = \left( 1 :0 : 0 : 0 : 1 \right), \\
 c_{7} & = \left( 1 : 0 : 0 : 1 : 0 \right), \\
c_{8} & =  \left( 1 :0 :0 : 0 : 0 \right). 
\end{align*}  

Let $ J = J_{0} \left( 42 \right) $ be the Jacobian of $X_{0} \left( 42 \right)$ and define 
\begin{align*}
D_{1} &  =  3 c_{1} -  c_{2} + 2c_{3} - 2c_{5} - 2 c_{6} ,   \\
D_{2}  &  =   -7 c_{1} +2 c_{2}  + 3c_{3}   +  c_{5} + c_{6} ,   \\
D_{3} &=   4c_{1} -4c_{2} + 3c_{3} - c_{4} -2c_{5}  - c_{6} + c_{8} ,    \\
D_{4} &=   6c_{1}  - 4 c_{2}  + 2c_{3}- 2 c_{5}  -2c_{6} .  \
\end{align*}

The linear equivalence classes of the above divisors are distinct, rational and we can verify using \texttt{Magma} that they are 2-torsion points on $J \cong \text{Pic}^{0} \left( X_{0} \left( 42 \right) \right)$.

Let $H$ be the subgroup of $J\left( K \right) \left[ 2 \right]$ generated by equivalence classes of elements of $HQT$ and the equivalence classes of $D_{1}, D_{2}, D_{3}, D_{4}$.
By factoring $11\mathcal{O}_{K}$ we find a prime ideal $\mathfrak{p} $ of norm 121. As 11 is a prime of good reduction for $X$, reduction modulo $\mathfrak{p}$ induces an injection 
\begin{center}
$r_{\mathfrak{p}} : J \left( K \right)_{\text{tors}} \longrightarrow J \left( \mathbb{F}_{\mathfrak{p}} \right)$. 
\end{center}
Using \texttt{Magma} we verify that  $J \left( \mathbb{F}_{\mathfrak{p}} \right) \cong \left( \mathbb{Z} / 2 \mathbb{Z} \right)^{4} \times \left( \mathbb{Z} / 8 \mathbb{Z} \right)^{2} \times \left( \mathbb{Z} / 24 \mathbb{Z} \right) \times \left( \mathbb{Z} / 48 \mathbb{Z} \right) \times \left( \mathbb{Z} / 192 \mathbb{Z} \right)^{2}$ and $r_{\mathfrak{p}} \left( H \right) \cong \left( \mathbb{Z} / 2 \mathbb{Z} \right)^{10}$, and thus the subgroup $H$ is necessarily the whole 2-torsion subgroup of $ J \left( K \right)_{\text{tors}} $.  Additionally, as the genus of $X$ is $5$, $H$ is in fact the entire 2-torsion subgroup 
\begin{center}
$H = J \left( \mathbb{\overline{Q}} \right)\left[ 2 \right].  $
\end{center} 
Let $G = \Gal \left( K / \mathbb{Q} \right) $. Elementary calculations show that $ G  \cong D_{8} \times C_{2} $, where $D_{8}$ is the dihedral group of order 8 and $C_{2}$ is the cyclic group of order 2.  Taking Galois invariants, we find 
\begin{center}
$H^{G} = J \left( \mathbb{Q} \right) \left[ 2 \right] \   =  \  \left( \mathbb{Z} / 2 \mathbb{Z} \right)\left[  D_{1} \right] +  \left( \mathbb{Z} / 2 \mathbb{Z} \right) \left[  D_{2} \right] +  \left( \mathbb{Z} / 2 \mathbb{Z} \right) \left[  D_{3} \right]  +  \left( \mathbb{Z} / 2 \mathbb{Z} \right) \left[ D_{4}  \right] $.
\end{center}

The rational cuspidal group of $J_0(42)$ was computed in \cite{OS} and found to be 
$$ C_0(42)(\Q) \cong \Z /2 \Z \times \Z / 2 \Z \times \Z / 12 \Z \times \Z / 48 \Z.$$
We further compute 
\begin{itemize}
    \item $J_0(42)(\mathbb{F}_5) \cong \Z /2 \Z \times \Z /4 \Z \times \Z /48 \Z \times \Z /48 \Z;$
    \item $J_0(42)(\mathbb{F}_{11}) \cong \Z / 2\Z \times \Z / 4 \Z \times \Z / 4 \Z \times \Z / 4 \Z \times \Z / 12 \Z \times \Z / 96 \Z.$
    \end{itemize}
By injectivity of torsion, $J_0(42)(\Q)_{\text{tors}}$ is isomorphic to either 
$$\Z /2 \Z \times \Z /2\Z \times \Z /12 \Z  \times \Z /48 \Z \ \text{or} \ \Z /2 \Z \times \Z / 4 \Z \times \Z / 12 \Z \times \Z / 48 \Z. $$
Suppose the latter holds. Then there exists $D \in J_0(42)(\Q)[4]$, which is not cuspidal. Since $J_0(42)(\Q)[2] = C_0(42)(\Q)[2]$, we necessarily have $2D \in C_0(42)(\Q)$. Working modulo $5$, we find $$ 2[D] = [3(c_5 -c_7)]. $$
However, no point of $J_0(42)(\mathbb{F}_{11})$ of order $4$ is such that twice it equals the reduction of $3(c_5-c_7)$, and thus by injectivity of torsion, no such $D$ exists, and thus $J_0(42)(\Q)_{\text{tors}} = C_0(42)(\Q).$
\section{Proof of Theorem $1$} 

To prove Theorem $1$ we begin by computing the two-torsion subgroup of $J_0(N)(\mathbb{Q})_{\text{tors}}$ for the five values of $N$ previously stated. The quadritangents planes required are stated in the \texttt{GitHub} directory given below. They were computed using the methods described in Sections \ref{Sch}, \ref{apr} and \ref{alg}. The points of their respective schemes of quadritangents were approximated using complex approximations and precise expressions were computed using the lattice technique described in Section \ref{CA}. For $N$ = $63$, $72$ and $75$ we were able to compute the entire $2$-torsion subgroup of the modular jacobian $J_{0}\left( N \right)$ using the computed quadritangents, and deduce the rational $2$-torsion subgroup as in the $N = 42$ case presented in detail in Section \ref{g5ex}. For details of these computations, the Galois orbits required in these computations and the \texttt{Magma} code used to compute them see 
\begin{center}
    \href{ https://github.com/ElviraLupoian/TwoTorsionSubgroups}{ https://github.com/ElviraLupoian/TwoTorsionSubgroups}
\end{center}
In some case we also used the cusps of the modular curves to simplify our computations. All the modular curves considered are non-hyperelliptic, of genus $5$ and their Jacobians have rank $0$ over $\Q$. The canonical models of the curves and cusps on these models used in our computations are as in \cite{OS}. 

In the case $N$ = $55$, we were unable to compute the entire $2$-torsion subgroup. We did compute sufficiently many quadritangents, defined over degree $12$ number field $K$, such that $J_{0}\left( 55 \right) \left( K \right) \left[ 2 \right] \cong \left( \Z / 2 \Z \right)^{8}$, and from this we deduced the rational $2$-torsion subgroup. The number field $K$ is not a Galois extension of $\Q$, and the computations are slightly different to those presented in \ref{g5ex}. An overview of this computation is given in Section \ref{X055}.

\subsection{$X_{0}\left( 55 \right)$}
\label{X055}
The canonical model of the curve $X = X_{0}\left( 55 \right)$ used to derive a scheme of quadritangents is the intersection of the three quadrics:
\begin{align*}
f_{1} & = x_{1}x_{3} - x_{2}^{2} + x_{2}x_{4} - x_{2}x_{5} - x_{3}^{2} + 3x_{3}x_{4} - 2x_{4}^{2} - 4x_{5}^{2}, \\
f_{2} & = x_{1}x_{4} - x_{2}x_{3} + 2x_{2}x_{4} - 2x_{2}x_{5} - 2x_{3}^{2} + 4x_{3}x_{4} + 5x_{3}x_{5}-2x_{4}^{2} - 4x_{4}x_{5} - 3x_{5}^{2}, \\
f_{3} & = x_{1}x_{5} - 2x_{2}x_{5} - x_{3}^{2} + 2x_{3}x_{4} + x_{3}x_{5} - x_{4}^{2}.
\end{align*}
and the cusps on this model are:
\begin{align*}
& c_{0} =  \left( -2 : 2 : 7 : 6 : 1 \right),\\
& c_{1} = \left( 3 : 2 : 2 : 1 : 1 \right), \\
& c_{2} =  \left( 1 : 0 :0 :0 :0 \right),\\  
& c_{3} = \left( 0 : 0 : 1 : 1 : 0 \right),
\end{align*}
as computed in \cite{OS}, where the authors show 
\begin{center}
 $C_{0}\left( 55 \right) \left( \Q \right) \cong \left( \Z / 10 \Z \right) \times \left( \Z / 20 \Z \right)$, \\
 $J_{0}\left( 55 \right) \left( \Q \right) / C_{0}\left( 55 \right) \left( \Q \right) \cong 0  $ or $ \left( \Z / 2 \Z \right)$ or $ \left( \Z / 2 \Z \right)^{2}$. 
\end{center}
Therefore to prove  $J_{0}\left( 55 \right) \left( \Q \right) =  C_{0}\left( 55 \right) \left( \Q \right)$, it is sufficient to show
\begin{center}
$J_{0}\left( 55 \right) \left( \Q \right) \left[ 2 \right]  =  C_{0}\left( 55 \right) \left( \Q \right) \left[ 2 \right] \cong \left( \Z/2 \Z \right)^{2}$.
\end{center}
Let $K$ be the number field defined by
\begin{center}
$f = u^{12} - 5u^{11} + 13u^{10} - 25u^9 + 39u^8 - 50u^7 + 53u^6 - 48u^5 + 37u^4 - 23u^3 + 12u^2 -4u + 1.$ 
\end{center}
Using the method described in Sections \ref{Sch3}, \ref{apr2} and \ref{CA} we find 18 quadritangent planes defined over $K$. The equations for these can be found in the \texttt{GitHub} directory stated at the beginning of this section.

Let $H$ be the subgroup of $2$-torsion points obtained from the divisors of the ratios of the equations defining these $18$ quadritangents. The ideal $\mathfrak{p} = \langle 47, 36 + \theta \rangle $, where $\theta$ is a generator of $K$, is prime of norm $47$, and the map induced by reduction modulo $\mathfrak{p}$
\begin{center}
$r_{\mathfrak{p}} : J_{0}\left( 55 \right) \left( K \right) \longrightarrow J_{0}\left( 55 \right) \left( \mathbb{F}_{\mathfrak{p}} \right) $
\end{center}
is an injection on the torsion subgroup of $ J_{0}\left( 55 \right) \left( K \right) $. Using \texttt{Magma}, we find 
\begin{center}
$J_{0}\left( 55 \right) \left( \mathbb{F}_{\mathfrak{p}} \right) \cong \left( \Z / 2 \Z \right)^{5} \times \left( \Z / 20 \Z \right)^{2} \times \left( \Z / 17220 \Z \right)$,
\end{center}
the image of $H$ is $r_{\mathfrak{p}} \cong \left( \Z / 2 \Z \right)^{8}$, and hence 
\begin{center}
 $  J_{0}\left( 55 \right) \left( K \right) \left[ 2 \right] = H \cong \left( \Z / 2 \Z \right)^{8}.$ 
\end{center}
Notably $K$ is not a Galois extension of $\mathbb{Q}$, as its automorphism group is a cyclic group of order $6$. The subfield of $K$ fixed by its automorphism group is $\Q \left( \sqrt{-11} \right)$, and thus $ \Q \left( \sqrt{-11} \right) \subset K $ is a Galois extension of degree $6$. Let $G = \Gal \left( \Q \left( \sqrt{-11} \right) / \Q \right)$. The $2$-torsion points of  $J_{0}\left( 55 \right) \left( \Q \left( \sqrt{-11} \right) \right)$ are simply the $2$-torsion points of $J_{0}\left( 55 \right) \left( K \right)$ fixed by $G$, and we find 
\begin{center}
 $ J_{0}\left( 55 \right) \left( \Q \left( \sqrt{-11} \right) \right) \left[ 2 \right] \cong \left( \Z / 2 \Z \right)^{3}$   
\end{center}

The action of $G$ on the 18 quadritangants has four orbits of size $6$, $6$, $3$ and $3$. We label these as $\{ P_{1}, \ldots, P_{6} \}$, $\{ Q_{1}, \ldots, Q_{6} \}$, $\{ R_{1}, R_{2}, R_{3} \}$, $\{ R_{4}, R_{5}, R_{6} \}$, where the elements of the above sets represent the defining equations of the quadritangents in each orbit.

Let 
\begin{align*}
F_{1} & = P_{1}\times \ldots \times P_{6}, \\
F_{2} & = Q_{1}\times \ldots \times Q_{6}, \\
F_{3} &= R_{1} \times R_{2} \times R_{3}, \\
F_{4} &= R_{4} \times R_{5} \times R_{6}.
\end{align*}
These have coefficients belonging to $\Q \left( \sqrt{-11} \right)$. We denote by $\overline{F_{i}}$ the conjugate of $F_{i}$ by the non-trivial element of $ \tilde{G} = \Gal \left( \mathbb{Q}\left( \sqrt{-11} \  \right) / \mathbb{Q} \right)$, and let 
\begin{center} $D_{i} = \left[ \frac{1}{2} \text{div} \left( \frac{F_{i+1}}{F_{1}} \right) \right] $ for $ i = 1, 2,3;$ \\
 $\overline{D_{j}} =  \left[ \frac{1}{2} \text{div}\left( \frac{\overline{F_{j}}}{F_{1}} \right) \right] $ for $j = 1,2, 3, 4$.
\end{center}
These yield divisor classes that are $2$-torsion and the Galois action on these is clear. Using the cusps stated above and the map $r_{\mathfrak{p}}$ we can find generators of the $2$-torsion part of the cuspidal subgroup 
\begin{center}
 $C_{0}\left( 55 \right) \left( \Q \right) \left[ 2 \right] = \left( \mathbb{Z} /2 \mathbb{Z} \right) C_{1} + \left( \mathbb{Z} /2 \mathbb{Z} \right) C_{2} $,
 \end{center}
where $C_{1} = \left[-5c_{2} + 5c_{3}\right] $ and $C_{2} = \left[ -10c_{1} + 10c_{2} \right] $.
Let $\Tilde{H}$ be the subgroup of  $J_{0}\left( 55 \right) \left( \Q \left( \sqrt{-11} \right) \right) \left[ 2 \right]$ generated by $C_{1}$, $C_{2}$, $D_{i}$ with $i = 1, 2, 3$ and $\Bar{D_{j}}$ with $j= 1, 2, 3, 4$. The map $r_{\mathfrak{p}} $ is an injection when restricted to $\Tilde{H}$, and we use it to show $\Tilde{H} \cong \left( \Z / 2 \Z \right)^{3}$, and hence $\Tilde{H} = J_{0}\left( 55 \right) \left( \Q \left( \sqrt{-11} \right) \right) \left[ 2 \right]$. Taking Galois invariants we obtain $ J_{0}\left( 55 \right) \left( \Q \right) \left[ 2 \right] \cong \left( \Z / 2 \Z \right)^{2}$, and hence 
 \begin{center}
 $J_{0}\left( 55  \right) \left(  \mathbb{Q}  \right) \left[ 2 \right]  = C_{0}\left( 55 \right) \left( \Q \right) \left[ 2 \right]\left( \mathbb{Z} /2 \mathbb{Z} \right)C_{1} + \left( \mathbb{Z} /2 \mathbb{Z} \right)C_{2} $.
 \end{center}
This completes our calculation. 
\begin{rem}
This is the only curve for which we were unable to calculate the whole $2$-torsion subgroup  $J_{0}\left( 55  \right) \left( \overline{ \mathbb{Q}} \right) \left[ 2 \right] \ \cong \ \left( \mathbb{Z} /2 \mathbb{Z} \right)^{10} $ . It is probable that the $2$-torsion subgroup is defined over the degree $24$ number field defined by
 \begin{center}
 $ f \left( u \right) = 888358082666609u^{24} - 686137237735072u^{23} + 965478109129036u^{22} - 
    108753611253152u^{21} + 1046333329183210u^{20} - 462274837855648u^{19} + 
    986194062109068u^{18} - 264174312478816u^{17} + 521023423262647u^{16} - 
    224217460467776u^{15} + 265604493047384u^{14} - 67790597560640u^{13} + 
    44563612667308u^{12} + 4777480913088u^{11} + 1939479463608u^{10} + 
    3337865504320u^{9} + 104055137263u^{8} + 362477031136u^{7} + 105446733532u^{6} +
    3486289120u^5 + 14677281802u^4 + 40615136u^3 + 680431932u^2 - 24640480u
    + 1394761
$
 \end{center}
 All the quadritangents found are defined over this number field. The field has few prime ideals of small norm making our computations impractical. 
 \end{rem}

\subsection{Excluding Small Torsion Orders}
Our computations show that $J_0(N)(\Q)[2] = C_0(N)(\Q)[2]$ for $N \in \{ 42, 55, 63, 72, 75\}$. To show that $C_0(N)(\Q) = J_0(N)(\Q)_{\text{tors}}$, we use our $2$-torsion calculation, alongside injectivity of torsion, see \cite{katz}. We note that the case of $N=42$ was treated in the previous section. 
\subsubsection{$N=55,63,75$} A simple reduction argument  suffices in these cases. For $N= 55$ we find that 
$$J_0(55)(\mathbb{F}_3) \cong \Z/2 \Z \times \Z / 2 \Z \times \Z /10 \Z \times \Z / 20 \Z.$$
As  $C_0(55)(\Q) \cong \Z /10 \Z \times \Z / 20 \Z$, by injectivity of torsion we find $J_0(55)(\Q)_{\text{tors}}$ is 
$$ C_0(55)(\Q) \ \text{or} \ \Z/ 2 \Z \times \Z / 10 \Z \times \Z / 20 \Z \ \text{or} \ \Z /2 \Z \times \Z / 2 \Z \times \Z / 10 \Z \times \Z / 20 \Z.$$
Our $2$-torsion calculation shows that $ J_0(55)(\Q)[2] \cong (\Z / 2 \Z)^{2},$ and thus the two latter possibilities cannot occur.  

For $N = 23$, reducing modulo $5$ and $23$ we find
\begin{itemize}
    \item $ J_0(63)(\mathbb{F}_5) \cong \Z / 2 \Z \times \Z / 2 \Z \times \Z / 2 \Z \times \Z / 2 \Z \times \Z / 4 \Z \times \Z / 96 \Z $;
    \item $J_0(63) (\mathbb{F}_{23}) \cong \Z / 6 \Z \times \Z / 12 \Z \times \Z / 48 \Z \times \Z / 2256 \Z$
\end{itemize}
and thus by injectivity of torsion $J_0(63)(\Q)_{\text{tors}}$ is $$C_0(63)(\Q) \cong \Z /2 \Z \times \Z / 4 \Z \times \Z / 48 \Z \ \text{or} \ \Z / 2 \Z \times \Z / 2 \Z  \times \Z / 4 \Z \times \Z / 48 \Z. $$
It follows from our $2$-torsion calculation that $J_0(63)(\Q)[2] = C_0(63)(\Q)[2] \cong (\Z / 2 \Z)^3$, and so the claim follows. 

For $N= 75$, we find:
\begin{itemize}
    \item $J_0(75)(\mathbb{F}_7) \cong \Z /2 \Z \times \Z / 4 \Z \times \Z / 8 \Z \times \Z /440 \Z;$
    \item $J_0(75)(\mathbb{F}_{11}) \cong \Z / 2 \Z \times \Z / 4 \Z \times \Z / 80 \Z \times \Z / 80 \Z. $
\end{itemize}
It follows that $J_0(75)(\Q)_{\text{tors}}$ is: 
$$ C_0(75)(\Q) \cong \Z / 2 \Z \times \Z / 4 \Z \times \Z / 40 \Z  \ \text{or} \  \Z / 2 \Z \times \Z / 4 \Z \times \Z / 4 \Z \times \Z / 40 \Z  \ \text{or} \Z / 2 \Z \times \Z / 4 \Z \times \Z / 8 \Z \times \Z / 40 \Z. $$
As our $2$-torsion calculation shows that $J_0(75)(\Q)[2] \cong (\Z / 2 \Z)^{3}$, the claim follows.

\subsubsection{$N= 72$} This case is treated in a similar manner to $N= 42$, and the calculation can be found in the \texttt{J072tor.m} file in the repository. 
 We work with the canonical model of the curve computed in \cite{OS}, that is, the genus $5$ curve in $\mathbb{P}^{4}$ cut out by the $3$ quadrics:
 \begin{itemize}
     \item[] $f_1 = x_1x_3 -x_2^2 + x_4^2; $
     \item[] $f_2 = x_1x_4 - x_3^2; $
     \item[] $ f_3 = x_1x_5 -x_3x_4-2x_5^2.$
 \end{itemize}
 By the work of Ozman-Siksek \cite{OS}, we find that $C_0(72)(\Q) \cong \Z / 2 \Z \times \Z / 4 \Z \times \Z / 12 \Z \times \Z / 12 \Z $ and $J_0(72)(\Q)_{\text{tors}}/ C_0(72)(\Q) \cong 0 \ \text{or} \ \Z/ 2 \Z$, and so $J_0(72)(\Q)_{\text{tors}}$ is isomorphic to one of the following groups:
 \begin{itemize}
     \item[1.] $\Z / 2 \Z \times \Z / 4 \Z \times \Z / 12 \Z \times \Z / 12 \Z;$
     \item[2.] $ \Z /4  \Z \times \Z / 4 \Z \times \Z /12 \Z \times \Z /12 \Z;$
     \item[3.] $\Z /2 \Z \times \Z /8 \Z \times \Z /12 \Z \times \Z /12 \Z;$
     \item[4.] $ \Z / 2 \Z \times \Z / 4 \Z \times \Z /12 \Z \times \Z /24 \Z;$
     \item[5.] $ \Z / 2 \Z \times \Z / 2 \Z \times \Z / 4 \Z \times \Z / 12 \Z \times \Z / 12 \Z.$ 
 \end{itemize}
 Our calculations show $J_0(72)(\Q)[2] \cong ( \Z / 2 \Z)^{4}$, which rules out (5) above. Moreover, we find that 

\begin{center}
    $ J_{0}(72)( \mathbb{F}_{5}) \cong \Z / 2 \Z \times \Z / 4 \Z \times \Z / 12 \Z \times \Z / 96 \Z  $
\end{center}
and hence, by the injectivity of torsion, $J_0(72)(\Q)_{\text{tors}}$ is congruent to either (1) or (4) above. 
We note that this further proves that $J_0(72)(\Q)[4] = C_0(72)(\Q)[4]$, and we find that linear equivalence classes of the following cuspidal divisors:
\begin{itemize}
    \item $D_1 = -5(2:0:0:0:1) + 2(4:3:2:1:1) + (4:-3:2:1:1) + (0:1:0:1:0) + (1:0:0:0:0);$
    \item $D_2 = 3(-2:0:2:-2:1) + 3(0:1:0:1:0) + 6(1:0:0:0:0) - 12 (2:0:0:0:1);$
    \item $D_3 = (1:0:-1:1:1) - (1:0:0:0:0)$
\end{itemize}
form a basis of this group. Suppose that $$ J_0(72)(\Q)_{\text{tors}} \cong \Z / 2 \Z \times \Z / 4 \Z \times \Z /12 \Z \times \Z /24 \Z.$$
Then, there exists $D \in J_0(72)(\Q)[8]$ which is not supported on cusps, but $2D \in C_0(72)(\Q)[4]$. Working modulo $5$, we find that such a $D$ necessarily satisfies:
$$2D = [ -D_1 + 2D_2 + D_3]. $$ 
Moreover, if such a $D$ exists, then for any prime $p \ge 7$, there exists $D' \in J_0(72)(\mathbb{F}_p)[8]$ such that $$2D'= [ -\tilde{D}_1 + 2 \tilde{D}_2 + \tilde{D}_3 ]$$
where $\tilde{a}$ denotes reduction modulo $p$. However, we find that this does not hold for any $8$-torsion point in $J_0(72)(\mathbb{F}_7)$. Thus 
$$ J_0(72)(\Q)_{\text{tors}} = C_0(72)(\Q) \cong \Z / 2 \Z \times \Z / 4 \Z \times \Z / 12 \Z \times \Z / 12 \Z.$$

 \bibliographystyle{abbrv}
\bibliography{ref}
\end{document}